\documentclass[12pt]{article}

\usepackage{latexsym}
\usepackage{amsfonts}
\usepackage{amsmath}
\usepackage{amssymb}
\usepackage{amscd}

\makeatletter
\def\@begintheorem#1#2{\it \trivlist
\item[\hskip \labelsep{\sc #1\ #2.}]}
\makeatother

\newcommand{\Hom}{\mathop{\rm Hom}\nolimits}
\newcommand{\ord}{\mathop{\rm ord}\nolimits}
\newcommand{\res}{\mathop{\rm res}\nolimits}
\newcommand{\im}{\mathop{\rm Im}\nolimits}

\newcommand{\ind}{\mathop{\rm ind}\nolimits}

\newcommand{\qed}{\nolinebreak\hfill\rule{2mm}{2mm}%
\par\medbreak}

\newcommand{\Supp}{\mathop{\rm Supp}\nolimits}
\newcommand{\proof}{\par\medbreak\noindent\it
Proof. \rm}

\newcommand{\FAT}[1]{\mbox{{$\mathbb{#1}$}}}

\newcommand{\cc}{\FAT{C}}
\newcommand{\qq}{\FAT{Q}}

\newcommand{\nn}{\FAT{N}}

\newtheorem{lem}{\sc Lemma}
\newtheorem{prop}[lem]{\sc Proposition}
\newtheorem{cor}[lem]{\sc Corollary}
\newtheorem{thm}[lem]{\sc Theorem}
\input{amssym.def}
\input{amssym.tex}
\title{Graded identities of matrix algebras and the universal
graded algebra}
\author{Eli Aljadeff \\
\small Department of Mathematics, Technion, Haifa 32000, Israel\\
\small aljadeff@tx.technion.ac.il\\
\\ 
Darrell Haile \\
\small Department of Mathematics, Indiana University, Bloomington,
IN 47405 \\
\small haile@indiana.edu\\
\\ 
Michael Natapov 
\\
\small Department of Mathematics, Indiana University, Bloomington,
IN 47405 \\
\small mnatapov@indiana.edu}

\date{}                                         

\begin{document}

\maketitle

\vspace{1cm}




\section{Introduction.}

In the last decade, group gradings and graded identities of finite
dimensional central simple algebras have been an active area of
research. We refer the reader to Bahturin, et al \cite{BSZ} and
\cite{BZ}. There are two basic kinds of group grading, elementary
and fine. It was proved by Bahturin and Zaicev \cite{BZ} that any
group grading of $M_{n}(\mathbb{C})$ is given by a certain
¨composition¨ of an elementary grading and a fine grading. In this
paper we are concerned with fine gradings on $M_{n}(\mathbb{C})$ and
their corresponding graded identities.

Let $R$ be a simple algebra, finite dimensional over its center $k$
and $G$ a finite group. We say that $R$ is fine graded by $G$ if
$R\cong \oplus_{g\in G} R_{g}$ is a grading and $\dim_{k}(R_{g})\leq
1$. Thus any component is either $0$ or isomorphic to $k$ as a
$k$--vector space. It is easy to show that $\Supp(R)$, the subset of
elements of $G$ for which $R_{g}$ is not $0$, is a subgroup of $G$.
Moreover $R$ is strongly graded by $\Supp(R)$, namely
$R_{g}R_{h}=R_{gh}$ for every $g,h \in \Supp(R)$. Since every group
$\Gamma$ containing $\Supp(R)$ provides a fine grading of $R$ (just
by putting $R_{g}=0$ for $g$ outside $\Supp(R)$) we will restrict
our attention to the case where $G=\Supp(R)$. In this case it has
been shown (see Bahturin and Zaicev \cite{BZ}) that the existence of
such a grading is equivalent to $R$ being isomorphic to a twisted
group algebra $k^cG$ where $[c]$ is an element in
$H^{2}(G,k^\times)$ and $G$ acts trivially on $k^\times$.

A group $G$ is said to be of central type if it admits a cocycle
$c$, $[c]\in H^{2}(G,\mathbb{C}^\times)$, for which the twisted
group algebra $\mathbb{C}^cG$ is central simple over $\mathbb{C}$.
By definition such a cocycle $c$ is called nondegenerate.  As we
have just seen fine gradings arise from such groups and cocycles.
Groups of central type appear in the theory of projective
representations of finite groups and in the classification of finite
dimensional Hopf algebras. Using the classification of finite simple
groups, Howlett and Isaacs \cite{HI} proved that any group of
central type is solvable.

Given a (fine) $G$--grading on $M_{n}(\mathbb{C})$ we consider the
set of graded identities: let $k$ be a subfield of $\mathbb{C}$ and
let $\Omega=\{x_{ig}:i\in \mathbb{N}, g\in G\}$ be a set of
indeterminates. Let $\Sigma(k)=k\langle \Omega \rangle$ be the
noncommutative free algebra generated by $\Omega$ over $k$. A
polynomial $p(x_{ig})\in \Sigma(k) $ is a graded identity of
$M_{n}(\mathbb{C})$ if the polynomial vanishes upon every
substitution of the indeterminates $x_{ig}$ by elements of degree
$g$ in $ M_{n}(\mathbb{C})$. It is clear that the set of graded
identities is an ideal of $\Sigma(k)$. Furthermore, it is a graded
$T$--ideal. Recall that an ideal (in $\Sigma(k)$) is a graded
$T$--ideal if it is closed under all $G$--graded endomorphisms of
$\Sigma(k)$. In section two we show that all graded identities are
already defined over a certain finite cyclotomic field extension
$\mathbb{Q}(\mu)$ of $\mathbb{Q}$. Moreover $\mathbb{Q}(\mu)$ is
minimal and unique with that property. We refer to $\mathbb{Q}(\mu)$
as the field of definition of the graded identities.  We consider
the free algebra $\Sigma(\mathbb{Q}(\mu))$ and denote by
$T(\mathbb{Q}(\mu))$ the $T$--ideal of graded identities of
$M_{n}(\mathbb{C})$ in $\Sigma(\mathbb{Q}(\mu))$. We introduce a set
of special graded identities which we call ``elementary''. We show
that a certain finite subset of the set of elementary identities
generate $T(\mathbb{Q}(\mu))$. In particular $T(\mathbb{Q}(\mu))$ is
finitely generated as a $T$--graded ideal.

In section three we examine the algebra $U_{G}=
\Sigma(\mathbb{Q}(\mu))/T(\mathbb{Q}(\mu))$, which we call the
universal $G$--graded algebra.  First we introduce another algebra
analogous to the ring of generic matrices in the classical theory.
We show this algebra is isomorphic to $U_G$ and thus are able to
prove that the center  $Z=Z(U_{G})$ is a domain and that if $F$
denotes the field of fractions of $Z$, then the algebra
$Q(U_G)=F\otimes_{Z}U_{G}$ is an $F$--central simple algebra of
dimension equal to the order of $G$. In particular $U_G$ is a prime
ring. Moreover we show there is a certain multiplicatively closed
subset $M$ in $Z$ such that the central localization $M^{-1}U_G$ is
an Azumaya algebra over its center $S=M^{-1}Z$.  The simple images
of this Azumaya algebra are the graded forms of $M_n(\mathbb{C})$,
that is the $G$--graded central simple $L$--algebras $B$ such that
$B\otimes_L\mathbb{C}$  is isomorphic as a graded algebra to
$M_n(\mathbb{C})$, where $L$ varies over all subfields of
$\mathbb{C}$. We also give quite explicit determinations of $S$ and
$F$.

Note that these algebras depend on the given grading on
$M_{n}(\mathbb{C})$ and so we should write $U_{G,c}$,
$M^{-1}U_{G,c}$ and $Q(U_{G,c})$, where $c$ is the given
nondegenerate two-cocycle.  We will omit the $c$ except in those
cases where we need to emphasize the particular grading.

Unlike the case of classical polynomial identities, the central
simple algebra $Q(U_{G})$  is not necessarily a division algebra.
Based on earlier work of the authors \cite{AHN} we in fact show that
$Q(U_{G,c})$ is a division algebra if and only if the group $G$
belongs to a certain explicit list $\Lambda$ of groups. This is
independent of the cocycle $c$, a fact we will return to below. The
list $\Lambda$ consists of a very special family of nilpotent
groups. Roughly speaking these are the nilpotent groups for which
each Sylow-$p$ subgroup is the direct product of an abelian group of
the form $A\times A$ (called of symmetric type) and possibly a
unique nonabelian group of the form $C_{p^{n}}\ltimes C_{p^{n}}$.
For $p=2$ an extra family of non abelian groups can occur, namely
$C_{2}\times C_{2^{n-1}}\ltimes C_{2^{n}}$. The precise definition
of $\Lambda$ is given in the last section.

The main result of the final section is that for every group $G$ on
the list $\Lambda$ the automorphism group of $G$ acts transitively
on the cohomology classes represented by nondegenerate two-cocycles.
It follows that the algebras $Q(U_{G,c})$  are all isomorphic for a
fixed $G$ (but not graded isomorphic).

For a group $G$ of central type let $\ind(G)$ denote the maximum
over all nondegenerate cocycles $c$ of the indices of the simple
algebras $Q(U_{G,c})$. We have just seen that if $G$ is not on the
list, then $\ind(G)$ is strictly less than the order of $G$.  In a
forthcoming paper, Aljadeff and Natapov \cite{AN}, it is shown that
in fact the groups on the list are the only groups responsible for
the index. The precise result is stated in the last section (Theorem
\ref{responsible.thm}). It follows from this theorem that if $G$ is
a $p$-group of central type then $\ind(G)\leq \max(\ord(H)^{1/2})$
where the maximum is taken over all $p$-groups $H$ on the list that
are sub-quotients of $G$.

If $Q(U_{G})$ is not a division algebra, there are graded identities
over the field of definition that are the product of nonidentities.
We present an explicit example of a grading on $M_6(\mathbb{C})$
(for the group $G=S_3\ltimes C_{6}$) for which there is a graded
identity over the field of definition that is the cube of a
nonidentity.

Even in the case where $Q(U_{G})$ is a division algebra it is possible
that this universal algebra does not remain a division algebra under
extension of the coefficient field.  In other words the algebra
$Q(U_{G})\otimes_{\mathbb{Q}(\mu)}\mathbb{C}$ may not be a division
algebra.   This means that there can be gradings on matrices for which
the algebra  $Q(U_{G})$ is a division algebra and yet there are graded
identities over $\mathbb{C}$
that are the products of nonidentities.  In
the last section we compute explicitly the index of
this extended algebra $Q(U_{G})\otimes_{\mathbb{Q}(\mu)}\mathbb{C}$
for groups $G$ on the list $\Lambda$.

\bigskip

\section{Graded Identities.}

Let $M_{n}(\mathbb{C})$ have a fine grading by the finite group $G$.
In this section we investigate the $T$--ideal of $G$--graded
identities on $M_{n}(\mathbb{C})$.  As we have seen in the
introduction the existence of such a grading implies that $G$ is of
central type and that there is a nondegenerate cocycle $c$ such that
$M_{n}(\mathbb{C})$ is isomorphic to the twisted group algebra
$\mathbb{C}^cG$.  We want to make this more precise.  Let
$A=M_{n}(\mathbb{C})$.  The $G$--grading induces a decomposition $A
\cong \oplus A_{g}$ where each homogeneous component is
$1$-dimensional over $\mathbb{C}$. Every homogeneous component
$A_{g}$ is spanned by an invertible element $u_{g}$ (which we also
fix from now on) and hence any element in $A_{g}$ is given by
$\lambda_{g}u_{g}$ where $\lambda_{g}\in \mathbb{C}$. The
multiplication is given by $u_{g}u_{h}=c(g,h)u_{gh}$ where
$c(g,h)\in \mathbb{C}^\times$ is a two-cocycle. We let $[c] \in
H^{2}(G,\mathbb{C}^\times)$ be the corresponding cohomology class.
In this way we identify $A$ (as a $G$-graded algebra) with the
twisted group algebra $\mathbb{C}^cG$. Replacing $c$ by a
cohomologous cocycle $c'$ produces a twisted group algebra
$\mathbb{C}^{c'}G$ that is isomorphic as a
$G$--graded algebra to $A$. 

We want to obtain information about the $T$--ideal of identities of
$M_{n}(\mathbb{C})$. We begin with a special family of identities.
As in the introduction we let $\Sigma(\mathbb{C})=\mathbb{C} \langle
\Omega \rangle$ denote the free algebra generated by $\Omega$ over
$\mathbb{C}$, where $\Omega=\{x_{ig}:i\in \mathbb{N}, g\in G\}$.


Let $Z_{1}=x_{r_1g_{1}}x_{r_2g_{2}}\cdots x_{r_kg_{k}}$ and
$Z_{2}=x_{s_1h_{1}}x_{s_2h_{2}}\cdots x_{s_lh_{l}}$ be two monomials
in $\Sigma(\mathbb{C})$.  We say $Z_1$ and $Z_2$ are {\it congruent}
if the following three properties are satisfied:

a) $g_{1} \cdots g_{k}=h_{1} \cdots h_{k}$

b) $k=l$ (equal length)

c) there exists an element $\pi \in Sym(k)$ such that for $1\leq i
\leq k$,  $ s_{i}=r_{\pi(i)}$ and $h_i=g_{\pi(i)}$.

We say $Z_{1}$ and $Z_{2}$ are {\it weakly congruent} if they
satisfy condition (a). Clearly congruence and weakly congruence are
equivalence relations.

Let $c(g_1,g_2,\dots, g_k)$ be the element in $\mathbb{C}^{\times}$
determined by $u_{g_1}u_{g_2}\cdots u_{g_k}=c(g_1,\dots,
g_k)u_{g_1\cdots g_k}$ in $\cc^cG$. For any two congruent monomials
$Z_{1}=x_{r_1g_{1}}x_{r_2g_{2}}\cdots x_{r_kg_{k}}$ and $Z_{2} = $
\\
$x_{r_{\pi(1)}g_{\pi(1)}}x_{r_{\pi(2)}g_{\pi(2)}}\cdots
x_{r_{\pi(k)}g_{\pi(k)}}$ consider the binomial

$$B = B(Z_1,Z_2) = Z_1-{c(g_1,\dots,g_k)\over
c(g_{\pi(1)},\dots,g_{\pi(k)})}Z_2$$

It is easy to see that $B$ is a graded identity. We will call it an
{\it elementary} identity. We refer to
$c(B)={c(g_1,\dots,g_k) \over c(g_{\pi(1)},\dots,g_{\pi(k)})}$ in
$\mathbb{C}^{\times}$ as the coefficient of the elementary identity
$B$.

This element $c(B)$ can be understood homologically.  Let
$F=F(\Omega)$ be the free group generated by $\Omega=\{x_{ig}:i\in
\mathbb{N}, g\in G\}$. Let $\nu: F\longrightarrow G$ $(x_{ig}
\longmapsto g)$ be the natural map and $R=\ker(\nu)$.  If $\pi\in
Sym(k)$ and $Z_1=x_{r_1g_{1}}x_{r_2g_{2}}\cdots x_{r_kg_{k}}$,
$Z_2=x_{r_{\pi(1)}g_{\pi(1)}}x_{r_{\pi(2)}g_{\pi(2)}}\cdots
x_{r_{\pi(k)}g_{\pi(k)}}$ are congruent, then the element
$z=Z_1(Z_2)^{-1}$ lies in  $[F,F]\cap R$.  By the Hopf formula the
Schur multiplier $M(G)$ of $G$ is given by  $M(G)=[F,F]\cap R/[F,R]$
and so we may consider the element $\overline{z}\in M(G)$. Moreover
the Universal Coefficients Theorem provides an isomorphism $\phi$
from
 $H^{2}(G,\mathbb{C}^\times)$ to $\Hom(M(G),\mathbb{C}^\times)$,

\begin{prop}
\label{elem.prop} Let the notation be as given above.

(a) The element $c(B)$ is independent of the representing cocycle $c$ of $[c]$.

(b) If $z=Z_1(Z_2)^{-1}$ is as given above then
$c(B)=\phi([c])(\overline{z})$.

(c) Each $c(B)$ is a root of unity, in fact $c(B)^{\ord(G)}=1$.

(d) If $J$ denotes the set of all elementary identities then the set
$\mu= \{c(B)\in \mathbb{C}^{\times}: B\in J\}$ is finite.
\end{prop}

\proof Statement (a)  is clear. For the proof of (b) we consider the
central extension

$$\{1\} \rightarrow R/[F,R] \rightarrow F/[F,R] \rightarrow G
\rightarrow \{1\}$$
induced by the map $\nu$ and let
$$\{1\} \rightarrow \mathbb{C}^\times \rightarrow \Gamma \rightarrow G
\rightarrow \{1\}$$ be the central extension that corresponds to the
cocycle $c$ (i.e. $\{u_{g}\}_{g\in G}$ is a set of representatives
in $\Gamma $ and $u_{g}u_{h}=c(g,h)u_{gh}$). It is easy to see that
the map $\gamma: F/[F,R] \rightarrow \Gamma $ given by
$\gamma(x_{ig})=u_{g}$ induces a map of extensions and it is well
known its restriction $\gamma: [F,F]\cap R/[F,R]=M(G)\rightarrow
\mathbb{C}^\times$ is independent of the presentation and
corresponds to $[c]$ by means of the Universal Coefficient Theorem.
That is $\phi([c])=\gamma$.  It follows that
$\phi([c])(\overline{z})=\gamma(\overline{z})= u_{g_{1}}\cdots
u_{g_{r}}(u_{g_{\pi(1)}}\cdots u_{g_{\pi(r)}})^{-1}$
$={c(g_1,\dots,g_k) \over c(g_{\pi(1)},\dots,g_{\pi(k)})}= c(B)$ as
desired.

Statements (c) and (d) are direct consequences of (b). This
completes the proof of the proposition.\qed

From the Proposition it follows that the set $J$ is defined over
$Q(\mu)$. Our next task is to show that the the $T$--ideal of graded
identities $T(\mathbb{C})$ is spanned by $J$ and that $Q(\mu)$ is
the field of definition for the graded identities. We will need the
following terminology:

Let $k$ be a subfield of  $\mathbb{C}$ and let $p(x_{r_ig_i})$$=$
$\sum \lambda_{r_1g_1,r_2g_2,\dots,r_kg_k}
x_{r_1g_1}x_{r_2g_2}\cdots x_{r_kg_k}$ be a polynomial in
$\Sigma(k)$ ($\subseteq$ $\Sigma(\mathbb{C})$). We say
$p(x_{r_ig_i})$ is {\it reduced} if the monomials
$x_{r_1g_1}x_{r_2g_2}\cdots x_{r_kg_k}$ are all different. If
$p(x_{r_ig_i})$ is reduced we say it is weakly homogeneous if its
monomials are all weakly congruent and we say $p(x_{r_ig_i})$ is
homogeneous is if its monomials are all congruent.  It is clear that
any reduced polynomial can be written uniquely as the sum of its
weakly homogeneous components and that each weakly homogeneous
component can be written uniquely as the sum of its homogeneous
components.

\begin{prop}
\label{decomposition.prop} If $p(x_{r_ig_i})$ is a graded identity
then its homogeneous components are graded identities as well.
Moreover, any graded identity is a linear combination of elementary
identities.
\end{prop}

\proof Let us show first that its weakly homogeneous components are
graded identities. Indeed, the replacement of $x_{r_ig_i}$ by
$u_{g_i}\in \cc^{c}G$ maps weakly homogeneous components of
$p(x_{r_ig_i})$ to different homogeneous components (in the graded
decomposition of $\cc^{c}G$). Since these are linearly independent,
they all must be $0$. Next, assume that we have at least two
homogeneous components in a weakly homogeneous component which we
denote by $\Pi$. Then there exists an indeterminate $y=x_{r_jg_{j}}$
which appears with different multiplicities in (at least) two
different monomials. Let $s\geq 1$ be the maximal multiplicity. We
decompose $\Pi$ into (at most) $s+1$ components $U_s,\dots,U_0$
where $U_i$ consists of the monomials of $\Pi$ that contain $y$ with
multiplicity $i$. Of course, by induction, the result will follow if
we show that $U_{s}$ is a graded identity. Suppose not. Then there
is a substitution which does not annihilate $U_{s}$. Note that since
$\Pi$ is weakly homogeneous the image of all components are
multiples of $u_g$ in $\cc^{c}G$ and hence there is an evaluation
that maps $U_i$ to $\lambda_i u_g$ with $\lambda_s$ not zero. We may
multiply the evaluation for $y$ by a central indeterminate $z$. Then
we get a non zero polynomial in $z$, whose coefficients are
$\lambda_i$ and for any evaluation of $z$ we get zero. This is of
course impossible in a field of characteristic zero.

To complete the proof of the proposition it suffices to show that
every homogeneous identity is a linear combination of elementary
identities. But this is clear since any two monomials which are
congruent determine an elementary identity (and monomials are not
identities).\qed

\begin{prop}
\label{fielddef.prop} If $L$ is a subfield of $\cc$ such that
$\mathbb{Q}(\mu)\subseteq L \subseteq \mathbb{C}$ then the set $J$
spans $T(L)$ over $L$. Conversely, if $L\subseteq \mathbb{C}$ is a
field of definition for $T(\mathbb{C})$, that is
$T(\mathbb{C})=T(L)\otimes _{L}\mathbb{C}$, then $L$ contains
$\mathbb{Q}(\mu)$.
\end{prop}

\proof Fix a monomial $Z=x_{r_1g_1}x_{r_2g_2}\cdots x_{r_kg_k}$. Let
$Z=Z_1, Z_2 \dots ,Z_d$ be the distinct monomials that are congruent
to $Z$ and let $W$ be the $d$-dimensional $\mathbb{C}$--vector space
spanned by these monomials. Each $Z_i$, for $1\leq i < d$ determines
an elementary identity $Z_i-c_iZ_d$ and these identities form a
basis for the  subspace $Y$ of graded identities in $W$.  In
particular the dimension of $Y$ over $\mathbb{C}$ is $d-1$.
Furthermore if $\gamma=\{c_1,\dots,c_{d-1}\}$ then clearly $Y$ is
defined over $\mathbb{Q}(\gamma)$. But $\gamma\subseteq \mu$ and so
it follows that $T(\mathbb{C})$ is defined over $\mathbb{Q}(\mu)$
and hence over any field $L$ that contains $\mathbb{Q}(\mu)$.

For the converse let $W_L$ denote the $L$--span of the monomials
$Z_1,\dots, Z_d$ and let $Y_L=Y\cap W_L$. We claim that if $L$ is a
field of definition for $Y$ (that is, $Y = Y_L\otimes_{L}
\mathbb{C}$), then $L\supseteq \mathbb{Q}(\gamma)$. Because $Z$ is
arbitrary, it will follow from the claim that $L$ contains $\mu$. To
prove the claim, let $f_{1},\dots,f_{d-1}$ be a basis of $Y_L$ over
$L$. Express $f_{1},\dots,f_{d-1}$ using the monomials $Z_{i}$. The
coefficient matrix is a $(d-1)\times d$ matrix over $L$. This matrix
may be row reduced to the normal form $(I_{d-1}, C^{'})$ where
$C^{'}$ is a $(d-1)\times 1$ matrix over $L$.  The column vector
$C^{'}$ is uniquely determined. In other words there are unique
scalars $a_i$ in $L$,  for $1\leq i <d$, such that the vectors
$Z_i-a_iZ_d$ form a basis for $Y_L$ over $L$. On the other hand the
identities $Z_i-c_iZ_d$ form a basis of $Y_{L(\gamma)}$ over
$L(\gamma)$. It follows that $a_i=c_i\in L$ for all $i$. \qed

We now show that the $T$--ideal of graded identities is finitely
generated. Let $n=\ord(G)$. Consider the set $V=\{x_{ig}: 1\leq i
\leq n\ ,  g \in G\}$ of indeterminates and let $E$ be the set of
elementary identities of length $\leq n$ (that is where the
monomials are of length $\leq n$) and  such that its indeterminates
are elements of $V$. Clearly $E$ is a finite set.

\begin{thm}
\label{finitegeneration.thm}
The ideal of graded identities is generated as a $T$--ideal by $E$.
In particular the ideal is finitely generated as a $T$--ideal.
\end{thm}

\proof Denote by $I'$ the $T$--ideal generated by $E$. We will show
that $I'$ contains all graded identities. Clearly it is enough to
show that $I'$ contains all elementary identities. Let
$B=x_{r_1g_{1}}x_{r_2g_{2}}\cdots x_{r_kg_{k}} -
c(B)x_{r_{\pi(1)}g_{\pi(1)}}x_{r_{\pi(2)}g_{\pi(2)}}\cdots
x_{r_{\pi(k)}g_{\pi(k)}}$ be an elementary identity of length $k$.
Clearly, if $k\leq n$, then $B$ is in $I'$. We assume therefore that
$k>n$ and proceed by induction on $k$. The first observation is that
if there exist $i,t$, with $1\leq i,t\leq k-1$ $i=\pi(t)$ and
$i+1=\pi(t+1)$ then we can reduce the length of the word:  We let
$g=g_ig_{i+1}$ and replace $x_{r_ig_i}x_{r_{i+1}g_{i+1}}$ by
$x_{rg}$, where $1\leq r\leq n$ and $x_{rg}$ does not appear in $B$.
The resulting identity (one has to check that the coefficient is
right) has shorter length and so we may assume this identity is in
$I'$.  But we can obtain the longer one from the shorter one by the
substitution of $x_{r_ig_i}x_{r_{i+1}g_{i+1}}$ for $x_{rg}$, so we
are done in this case.

Next, note that since $k>n$, the pigeonhole principle applied to the
expressions $g_1 \cdots g_s$, $s=1,\dots,k$ shows that there are
integers $1\leq i<j\leq n$ such that $g_ig_{i+1}\cdots g_j=1$.
Observe also, that if $\gamma$ is a cyclic permutation of the
numbers $i,i+1,\dots, j$ then  we have $g_{\gamma(i)}\cdots
g_{\gamma(j)}=1$.  That means that if we replace the string
$x_{r_ig_i}x_{r_{i+1}g_{i+1}}\cdots x_{r_jg_j}$ by a cyclic
permutation of the variables in the first monomial and leave the
second monomial alone, we do not change our identity modulo $I'$.
Also for all $g\in G$, we have $gg_ig_{i+1}\cdots
g_j=g_ig_{i+1}\cdots g_jg$. This means that modulo $I'$ we can move
the string $x_{r_ig_i}x_{r_{i+1}g_{i+1}}\cdots x_{r_jg_j}$ anywhere
in the monomial. So, combining these two things we can cyclically
move the substring and move it anywhere in the first monomial. But
by doing these moves we can arrange it so that some two variables
that are next to each other in the second monomial are also next to
each other (in the same order) in the first monomial. This reduces
us to the first case.\qed

\bigskip

\section{The Universal Algebra.} In this section we determine the
basic structural properties of the universal $G$--graded algebra
$U_{G}$ and relate that algebra to the other algebras described in
the introduction. The results are analogous to results in the theory
of (non-graded) polynomial identities. We will see for example that
$U_{G}$ is a prime ring and that its algebra of central quotients,
the universal $G$--graded algebra, is central simple.

We start with the following lemma.

\begin{lem}
\label{formsareforms.lem} Let $c$ be a nondegenerate cocycle on $G$
and let $\mathbb{C}^{c}G$ be the twisted group algebra. Let $L$ be a
subfield of $\mathbb{C}$.

(a) Let $\beta:G\times G\rightarrow L^{\times}$ be a two-cocycle and
let $L^{\beta}G$ be the twisted group algebra. There is a
homomorphism $\eta: L^{\beta}G \rightarrow \mathbb{C}^{c}G$ over $L$
of $G$--graded algebras if and only if the cocycles $c$ and $\beta$
are cohomologous in $\mathbb{C}$. In particular $\im(\eta)$ is a
$G$--graded form of $\mathbb{C}^{c}G$.

(b) There exists a two cocycle $\beta:G\times G\rightarrow
L^{\times}$ cohomologous to $c$ over $\mathbb{C}$ if and only if $L$
contains $\mu$, the set of Proposition~\ref{elem.prop}.
\end{lem}

\proof (a) If $\eta: L^{\beta}G \rightarrow \mathbb{C}^{c}G$ is a
homomorphism of $G$-graded $L$--algebras then there are bases
$\{v_{\sigma}\}\subset L^{\beta}G$, $\{u_{\sigma}\}\subset
\mathbb{C}^{c}G$, and scalars $\{\lambda_{\sigma}\}\subset
\mathbb{C^\times}$ such that
$v_{\sigma}v_{\tau}=\beta(\sigma,\tau)v_{\sigma\tau}$,
$u_{\sigma}u_{\tau}=c(\sigma,\tau)u_{\sigma\tau}$ and
$\eta(v_{\sigma})=\lambda_{\sigma}u_{\sigma}$. It follows that
$\beta$ and $c$ are cohomologous over $\mathbb{C}$. A similar
calculation shows the other direction.

(b) We have seen in Proposition~\ref{elem.prop} that $\phi([c])=\mu$
where $\phi$ is the isomorphism between $H^{2}(G,\mathbb{C}^\times)$
and $\Hom(M(G),\mathbb{C}^\times)$. Hence if $\beta$ is cohomologous
to $c$ over $\mathbb{C}$ then $L$ must contain $\mu$ by part (d) of
Proposition~\ref{elem.prop}.

For the converse we may assume that $L=\mathbb{Q}(\mu)$. By the
naturality of the Universal Coefficient Theorem we have a
commutative diagram

\[
\begin{CD}
1  @>>> Ext^{1}(G_{ab},\mathbb{Q}(\mu))   @>>> H^{2}(G,\mathbb{Q}(\mu)^\times)  @>>> \Hom(M(G),\mathbb{Q}(\mu)^\times)  @>>> 1\\
@.                 @.                     @VV\nu V                               @VV\overline{\nu}V \\
@.                   @.                   H^{2}(G,\mathbb{C}^\times)     @>\phi>>    \Hom(M(G),\mathbb{C}^\times)      \\
\end{CD}
\]

Clearly $\phi([c])\in \Hom(M(G),\mathbb{C}^\times)$ is in
$\im(\bar{\nu})$ and hence $[c]$ is in $\im(\nu)$.\qed

The first step in the analysis of $U_{G}$ is the construction of a
counterpart to the ring of generic matrices. Let $G$ be a group of
central type of order $n^2$ and let $c:G\times G\rightarrow
\mathbb{C}^\times$ be a nondegenerate two-cocycle. For every $g\in
G$ let $t_{ig}$ for $i=1,2,3,\dots$ be indeterminates. For each
$g\in G$ let $t_g=t_{1g}$. We will assume $u_1=1$ (i.e. the cocycle
$c$ is normalized). Let $k$ denote the field generated over $\Bbb
{Q}$ by the indeterminates and the values of the cocycle. Note that
by Lemma~\ref{formsareforms.lem}, $k$ contains $\mathbb{Q}(\mu)$.
Consider the twisted group algebra $k^{c}G$. Because the cocycle is
nondegenerate $\mathbb{Q}(c(g,h))^{c}G$ is a central simple
$\mathbb{Q}(c(g,h))$--algebra and hence $k^{c}G$ is a central simple
$k$--algebra.  Now let  $\overline{U_G}$ denote the
$\Bbb{Q}(\mu)$--subalgebra of $k^cG$ generated by the elements
$t_{ig}u_{g}$ for all $i$ and all $g\in G$.  We want to describe
$\overline{U_G}$.  Let $\overline{Z}$ denote the center of
$\overline{U_G}$.

\begin{prop}
\label{centerofgeneric.prop} (a) The center  $\overline{Z}$ of
$\overline{U_G}$ is the $\mathbb{Q}(\mu)$--subalgebra of
$\overline{U_G}$ generated by the following set of elements:

$$\{t_{i_1g_1}u_{g_1}t_{i_{2}g_2}u_{g_2}\cdots t_{i_{m}g_m}u_{g_m}
\ | \  g_1g_2\cdots g_m=1 \}$$

(b) We have  $\overline{Z}\subseteq k$.  In particular  $\overline{Z}$ is an integral
domain.
\end{prop}

\proof The ring $k\overline{U_G}$ contains the elements $u_g$ for
all $g\in G$ and so $k\overline{U_G}=k^cG$, which has center $k$.
Hence the center of $\overline{U_G}$ is $\overline{U_G}\cap k$. But
this intersection is precisely the $\mathbb{Q}(\mu)$ span of the
given set of monomials. This proves both parts.\qed

Now let $s:G\times G\rightarrow k^\times$ be the function given by
$s(g,h)={t_gt_h \over t_{gh}}c(g,h)$.  This is a two-cocycle,
cohomologous to $c$ over $k$.  Because $c$ is nondegenerate, so is
$s$ and hence the algebra $k^sG$ is $k$--central simple. In fact it
is convenient to view the algebra $k^sG$ as equal to $k^cG$. It is
spanned over $k$ by the elements $t_{g}u_{g}$.

Let $Y$ be the subgroup of $k^\times$ generated by the values of $s$
and let $\overline{S}=\mathbb{Q}(\mu)[Y]$, a subring of $k$. The
$\overline{S}$--subalgebra of $k^sG$ generated by the elements
$t_{g}u_{g}$  is the twisted group algebra $\overline{S}^sG$.  Let
$\overline{M}=\{t_{g_1}u_{g_1}t_{g_2}u_{g_2}\cdots t_{g_m}u_{g_m} \
| \ g_1g_2\cdots g_m=1 \}$.  Let $\overline{F}\subseteq k$ denote
the field of fractions of $\overline{Z}$.

\begin{prop}
\label{localizationofcenter.prop} (a) The set $\overline{M}$ is a
multiplicatively closed subset of $\overline{Z}$ and
$\overline{M}^{-1}\overline{Z}=\overline{S}[{t_{ig} \over t_g} \ | \
i\geq 1, g\in G]$.

(b) We have
$\overline{U_G}(\overline{M}^{-1}\overline{Z})=\overline{S}^sG(\overline{M}^{-1}\overline{Z})$.
Moreover
$\overline{U_G}\overline{F}=\overline{S}^sG\overline{F}=\overline{F}^sG$.

(c) $\overline{F}^sG$ is a central simple $\overline{F}$--algebra.

(d) The ring $\overline{U_G}$ is a prime ring with ring of central
quotients $Q(\overline{U_G})$ isomorphic to the central simple
algebra $\overline{F}^sG$.
\end{prop}

\proof (a)  We first show $\overline{M}^{-1}\overline{Z}\supseteq
\overline{S}[{t_{ig} \over t_g} \ | \ i\geq 1, g\in G]$. Observe
that $t_{e}$ and $t_{g}t_{g^{-1}}c(g,g^{-1})$ are in $\overline{M}$
and hence
$t_{ig}t_{g^{-1}}c(g,g^{-1})/t_{g}t_{g^{-1}}c(g,g^{-1})=t_{ig}/t_{g}$
is in $\overline{M}^{-1}\overline{Z}$. Next, the elements
$s(g,h)^{\pm 1}$ are in $\overline{M}^{-1}\overline{Z}$ because the
elements

$t_{g}u_{g}t_{h}u_{h}t_{(gh)^{-1}}u_{(gh)^{-1}}=t_{g}t_{h}t_{(gh)^{-1}}c(g,h)c(gh,(gh)^{-1})$
and

$t_{gh}u_{gh}t_{(gh)^{-1}}u_{(gh)^{-1}}=t_{gh}t_{(gh)^{-1}}c(gh,(gh)^{-1})$
are in $\overline{M}$.  Similar calculations show that the opposite inclusion
also holds.

(b) This follows easily from part (a).

(c) We have $\overline{F}^sG\otimes_{\overline{F}}k=k^sG$ which is
$k$--central simple and so $\overline{F}^sG$ is
$\overline{F}$--central simple.

(d) By part (b) we have
$\overline{U_G}\overline{F}=\overline{S}^sG\overline{F}=\overline{F}^sG$
and $\overline{F}^sG$ is $\overline{F}$--central simple, the ring
$\overline{U_G}$ is prime and because $\overline{F}$ is the ring of
fractions of $\overline{Z}$ we see that the ring of central
quotients $Q(\overline{U_G})$ is isomorphic to $\overline{F}^sG$.
\qed

Let $\psi$ denote the $\mathbb{Q}(\mu)$--algebra homomorphism from
the free algebra $\Sigma(\mathbb{Q}(\mu))$ to $k^{c}G$ given by
$\psi(x_{ig})=t_{ig}u_{g}$. The image of this map is clearly
$\overline{U_G}$.

\begin{prop}
\label{iso.prop} The kernel of $\psi$ is the ideal
$T(\mathbb{Q}(\mu))$ of graded identities of $\mathbb{C}^{c}G$ and
hence $\psi$ induces a $\mathbb{Q}(\mu)$--algebra isomorphism from
$U_G$ to $\overline{U_G}$. In particular the universal algebra
$U_G$ is a prime ring with center isomorphic to $\overline{Z}$ and and its
ring of central quotients, that is $Q(U_G)$ is isomorphic
to $\overline{F}^sG$.
\end{prop}

\proof This is clear since a polynomial $p(x_{ig})\in
\Sigma(\mathbb{Q}(\mu))$ is an identity of $\mathbb{C}^{c}G$ if
and only if $p(\lambda_{ig}u_{g})=0$ for any $\lambda_{ig}\in
\mathbb{C}$ and this is equivalent to $p(t_{ig}u_{g})=0$ where
$\{t_{ig}\}$ are central indeterminates. \qed

Because of this isomorphism we will from now on drop the
bars on $\overline{Z}$, $\overline{F}$, etc.

We want to say more about $F$, the field of fractions of the
center of the universal algebra $U_G$.  The group $Y$ is finitely
generated and hence of the form $Y_tY_f$ where $Y_t$ is the
torsion subgroup of $Y$ and $Y_f$ is a finitely generated free
abelian group.  Because $Y$ is a subgroup of $k^\times$ it follows
that $Y_t$ is cyclic. We will see in the next proposition that in
fact $Y_t=\mu$.  Let $y_1,y_2,\dots ,y_m\in k$ denote a basis for
$Y_f$. (We will see later - in Proposition~\ref{rank.prop}  - that the
rank of $Y_f$ is equal to $n$, the order of $G$, and so $m=n$.)

\begin{prop}
\label{fieldF.prop} (a) We have $Y_t=Y\cap \mathbb{Q}(\{c(g,h) \ | \
g,h\in G\})=\mu$.

(b) We have $S=\mathbb{Q}(\mu)[Y]=\mathbb{Q}(\mu)[y_1^{\pm
1},y_2^{\pm 1},\dots, y_m^{\pm 1}]$, the ring of Laurent polynomials
in $y_1,\dots, y_m$ over $\mathbb{Q}(\mu)$.

(c) We have $F=\mathbb{Q}(\mu)(y_1,\dots, y_m)({t_{ig} \over t_g} \
| \  i\geq 1, g\in G)$.  In particular $F$ is isomorphic to the
field of rational functions in countably many variables over
$\mathbb{Q}(\mu)$.
\end{prop}

\proof To prove (a) we prove the following inclusions: (i) $Y_t
\subseteq Y\cap \mathbb{Q}(\{c(g,h) \ | \ g,h\in G\})$, (ii) $\mu
\subseteq Y_t$, and (iii) $Y\cap \mathbb{Q}(\{c(g,h) \ | \ g,h\in
G\}) \subseteq \mu$.  To prove (i) consider an arbitrary element $z$
in $Y$, $z=\prod_{i=1}^{n}
s(g_{i},h_{i})^{\epsilon_{i}}=\prod_{i=1}^{n} ({t_{g_{i}}t_{h_i}
\over t_{g_ih_i}}c(g_i,h_i))^{\epsilon_{i}}$ where ${\epsilon_{i}}$
are $\pm 1$. If the product is not independent of the $t's$ (say
$t_{g}$ appears) it is clear that any positive power of $z$ contains
a power of $t_{g}$ and hence does not equal 1. This shows (i). Next,
let $\lambda\in \mu$ and let
$$x_{1g_{1}}x_{1g_{2}}\cdots x_{1g_{k}}- \lambda
x_{1g_{\pi(1)}}x_{1g_{\pi(2)}}\cdots x_{1g_{\pi(k)}}$$

\noindent be an elementary identity. It follows that
$$t_{g_{1}}u_{g_{1}}t_{g_{2}}u_{g_{2}}\cdots t_{g_{k}}u_{g_{k}}= \lambda
t_{g_{\pi(1)}}u_{g_{\pi(1)}}t_{g_{\pi(2)}}u_{g_{\pi(2)}}\cdots
t_{g_{\pi(k)}}u_{g_{\pi(k)}}$$ and hence
$s(g_{1},g_{2})s(g_{1}g_{2},g_{3})\cdots s(g_{1}g_{2}\cdots
g_{k-1}, g_{k})t_{g_{1}g_{2}\cdots g_{k}}u_{g_{1}g_{2}\cdots
g_{k}}=$

\noindent $\lambda
s(g_{\pi(1)},g_{\pi(2)})s(g_{\pi(1)}g_{\pi(2)},g_{\pi(3)})\cdots
s(g_{\pi(1)}g_{\pi(1)}\cdots g_{\pi(k-1)},
g_{\pi(k)})t_{g_{\pi(1)}g_{\pi(2)}\cdots
g_{\pi(k)}}u_{g_{\pi(1)}g_{\pi(2)}\cdots g_{\pi(k)}}$.

\noindent But $g_{1}g_{2}\cdots g_{k}=g_{\pi(1)}g_{\pi(2)}\cdots
g_{\pi(k)}$ in $G$ and so $\lambda$ lies in $Y$. Because $\lambda$
is a root of unity, we get $\lambda\in Y_t$. This proves (ii).

To prove (iii) let $\lambda$ be an element of $Y\cap
\mathbb{Q}(\{c(g,h) \ | \ g,h\in G\})$. We can write $\lambda=
s(g_{1},h_{1})^{\epsilon_{1}}s(g_{2},h_{2})^{\epsilon_{2}} \cdots
s(g_{k},h_{k})^{\epsilon_{k}}=$

$({t_{g_{1}}t_{h_{1}} \over
t_{g_{1}h_{1}}}c(g_{1},h_{1}))^{\epsilon_{1}}({t_{g_{2}}t_{h_{2}}
\over t_{g_{2}h_{2}}}c(g_{2},h_{2}))^{\epsilon_{2}} \cdots
({t_{g_{k}}t_{h_{k}} \over
t_{g_{k}h_{k}}}c(g_{k},h_{k}))^{\epsilon_{k}}$.

\noindent We use the notation and set-up of part (b) of
Proposition~\ref{elem.prop}.  Let $d$ be the element of the free
group $F$ given by
$d=x_{g_{1}}x_{h_{1}}(x_{g_{1}h_{1}})^{-1}x_{g_{2}}x_{h_{2}}(x_{g_{2}h_{2}})^{-1}\cdots
x_{g_{k}}x_{h_{k}}(x_{g_{k}h_{k}})^{-1}$.  Clearly $d$ lies in $R$.
Because $\lambda$ lies in $\mathbb{Q}(\{c(g,h) \ | \ g,h\in G\})$,
the variables $t_g$ in the expression for $\lambda$ must all cancel
and so $d\in [F,F]$ and
$\lambda=c(g_{1},h_{1}))^{\epsilon_{1}}c(g_{2},h_{2}))^{\epsilon_{2}}
\cdots c(g_{k},h_{k}))^{\epsilon_{k}}$. In particular we may
consider the element $\overline{d}$ in $M(G)$.  Then again as in
part (b) of Proposition~\ref{elem.prop} we have
$\phi([c])(\overline{d})=\gamma(\overline{d})=c(g_{1},h_{1}))^{\epsilon_{1}}c(g_{2},h_{2}))^{\epsilon_{2}}
\cdots c(g_{k},h_{k}))^{\epsilon_{k}}=\lambda$ and so $\lambda$ lies
in $\mu$.

Part (b) is clear. Part (c) follows from part (b) and from the fact
that $M^{-1}Z=S[{t_{ig} \over t_g} \ | \ i\geq 1, g\in G]$.\qed

Our next goal is to study the algebra $S^sG$.  We will see that it
is Azumaya and that its simple images are the graded forms of
$\mathbb{C}^cG=M_n(\mathbb{C})$.

\begin{prop}
\label{Az.prop} (a) Let $m$ be a maximal ideal of $S$ and let
$y\rightarrow \overline{y}$ denote the canonical homomorphism from
$S$ to $\tilde{S}=S/m$. Consider the two-cocycle $\tilde{s}:G\times
G \longrightarrow S/m$ induced by $s$ and let
$\tilde{S}^{\tilde{s}}G$ be the corresponding twisted group algebra.
Then the algebra $S^sG/mS^sG$ is isomorphic to
$\tilde{S}^{\tilde{s}}G$ and is a central simple
$\tilde{S}$--algebra.

(b)  The ring $S^sG$ is Azumaya over $S$.
\end{prop}

\proof (a) We first show that the two-cocycle $\tilde{s}$ is
nondegenerate.  We use the following characterization of
nondegeneracy (see Isaacs \cite[Problem 11.8]{I}): a two-cocycle
$\delta:G\times G\rightarrow \mathbb{C}^\times$ is nondegenerate if
and only if for every element $g\in G$ there is an element $h\in G$
such that $g$ and $h$ commute but in the algebra $\mathbb{C}^\delta
G$ the elements $u_g$ and $u_h$ do not commute. Let $g\in G$.
Because $c$ is nondegenerate there is an element $h\in G$ that
centralizes $g$ and such that $c(g,h)c(h,g)^{-1}\neq 1$. Moreover
$c(g,h)c(h,g)^{-1}\in \mu$ and $s(g,h)s(h,g)^{-1}=
c(g,h)c(h,g)^{-1}$. Because $S$ contains the field
$\mathbb{Q}(\mu)$, the restriction of the canonical homomorphism
from $S$ to $\tilde{S}$ is injective on $\mu$ and so
$\tilde{s}(g,h)\tilde{s}(h,g)^{-1}\neq 1$.  Hence $\tilde{s}$ is
nondegenerate.  It follows that $\tilde{S}^{\tilde{s}}G$ is a
central simple $\tilde{S}$--algebra.  It is clear that the canonical
homomorphism from $S^sG$ to  $\tilde{S}^{\tilde{s}}G$ is surjective
with kernel $mS^sG$.

(b) This follows from part (a) (see DeMeyer and Ingraham
\cite[Theorem 7.1]{DI}).\qed

Let $L$ be any subfield of $\mathbb{C}$.  Recall that every
$G$--graded form of $\mathbb{C}^cG$  over $L$ is a twisted group
algebra $L^\beta G$ where $\beta$ is a two-cocycle cohomologous to
$c$ over $\mathbb{C}$.

We will show now that $S^{s}G$ is universal with respect to such
forms.  The following observation about the Hopf formula
$M(G)=[F,F]\cap R/[F,R]$ will be useful (and is probably well
known). We use the notation and discussion preceding
Proposition~\ref{elem.prop}.

\begin{lem}
\label{Hopf.lem}
Every element in $[F,F]\cap R$ is equivalent modulo $[F,R]$ to an
element of the form $x_{g_{1}} \cdots x_{g_{k}}(x_{g_{\pi (1)}}
\cdots x_{g_{\pi (k)}})^{-1}$ where  $\pi\in Sym(k)$ and  $g_{1}\cdots
g_{k}=g_{\pi (1)} \cdots g_{\pi (k)}$.
\end{lem}

\proof  Let $y\in [F,F]\cap R$.  We can write
$y=x_{g_1}^{\epsilon_1}x_{g_2}^{\epsilon_2}\cdots
x_{g_m}^{\epsilon_m}$ where $\epsilon_i=\pm 1$ for $1\leq i\leq m$
and where for each group element $g$ that appears the elements $x_g$
and $x_g^{-1}$ appear the same number of times (and so $m$ is even)
and where the product $g_1^{\epsilon_1}g_2^{\epsilon_2}\cdots
g_m^{\epsilon_m}$ equals one. Let $m=2r$.  We \underline{claim} that
we may assume $\epsilon_i=1$ if $1\leq i\leq r$ and $\epsilon_i=-1$
if $r+1\leq i\leq m$.  The lemma follows easily from this claim.  To
prove the claim, suppose $\epsilon_i=-1$ for some $i$, $1\leq i\leq
r$.  We consider $y=x_{g_1}^{\epsilon_1}x_{g_2}^{\epsilon_2}\cdots
x_{g_m}^{\epsilon_m}=x_{g_1}^{\epsilon_1}x_{g_2}^{\epsilon_2}\cdots
x_{g_r}^{\epsilon_{r}}(x_{g_i}x_{g_i^{-1}}(x_{g_i^{-1}})^{-1}(x_{g_i})^{-1})x_{g_{r+1}}^{\epsilon_{r+1}}\cdots
x_{g_m}^{\epsilon_m}$.  Now for any element $g\in G$ the element
$x_gx_{g^{-1}}$ commutes modulo $[F,R]$ with every element of $F$.
Hence we can move the element $x_{g_i}x_{g_i^{-1}}$ in $y$ and
obtain an element that is equivalent modulo $[F,R]$ to $y$ and for
which the element $x_{g_i}^{-1}$  is replaced by $x_{{g_i}^{-1}}$.
Continuing in this way proves the claim.  \qed

\begin{prop}
\label{forms.prop} Let $L$ be a subfield of $\mathbb{C}$ which
contains $\mathbb{Q}(\mu)$. Let $\varphi:S\longrightarrow L$ be a
$\mathbb{Q}(\mu)$--algebra homomorphism.  Let $\beta(g,h)=\varphi
(s(g,h))$. Then $\varphi$ induces a $G$--graded homomorphism
$\hat{\varphi}: S^{s}G\longrightarrow L^{\beta}G$ and $L^\beta G$ is
a $G$--graded form of $\mathbb{C}^cG$. In particular $\beta$ is
cohomologous to $c$ over $\mathbb{C}$.

Conversely if $L^\beta G$ is a $G$--graded form of $\mathbb{C}^cG$,
then there is homomorphism $\varphi$ of $S$ into $L$ such that
$\gamma=\varphi (s(g,h)$ is a cocycle on $G$ cohomologous to $\beta$
over $L$ and $\varphi$ induces a homomorphism from $S^sG$ onto
$\varphi(S)^\gamma G$, a $G$--graded form of $L^\beta G$.
\end{prop}

\proof Let  $\varphi:S\longrightarrow L$ be a
$\mathbb{Q}(\mu)$--algebra homomorphism and let $\beta(g,h)=\varphi
(s(g,h))$.  We have seen in Proposition~\ref{Az.prop} that $\beta$
is a nondegenerate cocycle on $G$ and hence that $L^\beta G$ is a
central simple $L$--algebra.  Clearly $\varphi$ induces a
homomorphism $\hat{\varphi}: S^{s}G\longrightarrow L^{\beta}G$ and
if $A$ denotes the image of $S^s G$ and $A$ has center $R$ (say),
then $A\otimes_RL=L^\beta G$. We are left with showing $L^\beta G$
is a form of $\mathbb{C}^cG$, that is we need to show  $\beta$ and
$c$ are cohomologous over $\mathbb{C}$.  It suffices to show that
they determine the same element in $\Hom(M(G),\mathbb{C}^\times)$,
that is that $\phi([c])=\phi([\beta])$, where $\phi$ is the
isomorphism from $H^{2}(G,\mathbb{C}^\times)$ to
$\Hom(M(G),\mathbb{C}^\times)$ (see the discussion preceding
Proposition~\ref{elem.prop}). By the lemma it suffices to show that
if $g_{1}\cdots g_{k}=g_{\pi (1)} \cdots g_{\pi (k)}$ then
$\phi([\beta])(y)=\phi([c])(y)$ where $y=(x_{g_{1}} \cdots
x_{g_{k}})(x_{g_{\pi (1)}} \cdots x_{g_{\pi (k)}})^{-1}$.  For the
purposes of the proof it will be useful to use the following
notation: If $\alpha$ is any two-cocycle then we denote the
expression
$${\alpha(g_{1},g_{2})\alpha(g_{1}g_{2},g_{3}) \cdots
\alpha(g_{1}g_{2} \dots g_{k-1},g_k) \over \alpha(g_{\pi (1)},g_{\pi
(2)})\alpha(g_{\pi (1)}g_{\pi (2)},g_{\pi (3)}) \cdots \alpha(g_{\pi
(1)}g_{\pi (2)} \dots g_{\pi (k-1)},g_{\pi (k)})}$$ by
$\alpha(g_1,g_2,\dots ,g_k)/\alpha(g_{\pi (1)},g_{\pi (2)}, \dots,
g_{\pi (k)})$.

We have $\phi([\beta])(y)=\beta(g_1,g_2,\dots ,g_k)/\beta(g_{\pi
(1)},g_{\pi (2)}, \dots, g_{\pi (k)})=$

\noindent $\varphi (s(g_1,g_2,\dots ,g_k)/s(g_{\pi (1)},g_{\pi (2)},
\dots, g_{\pi (k)}))=$

\noindent $\varphi(
(t_{g_1}t_{g_2}/t_{g_1g_2})(t_{g_1g_2}t_{g_3}/t_{g_1g_2g_3})\cdots (t_{g_1g_2\cdots
    g_{k-1}}t_{g_k}/t_{g_1g_2\cdots g_k})/$

$(t_{g_{\pi (1)}}t_{g_{\pi (2)}}/t_{g_{\pi (1)}g_{\pi (2)}})\cdots
(t_{g_{\pi (1)}g_{\pi (2)}\cdots g_{\pi (k-1)}}t_{g_{\pi (k)}}/
t_{g_{\pi (1)}g_{\pi (2)}\cdots g_{\pi (k)}})$

$c(g_1,g_2,\dots ,g_k)/c(g_{\pi (1)},g_{\pi (2)}, \dots, g_{\pi
(k)}))=$

\noindent $\varphi( (t_{g_1g_2\cdots g_k}/t_{g_{\pi (1)}g_{\pi
(2)}\cdots g_{\pi (k)}}) c(g_1,g_2,\dots ,g_k)/c(g_{\pi (1)},g_{\pi
(2)}, \dots, g_{\pi (k)}))=$

\noindent $\varphi(c(g_1,g_2,\dots ,g_k)/c(g_{\pi (1)},g_{\pi (2)},
\dots, g_{\pi (k)}))$, because $g_{1}\cdots g_{k}=g_{\pi (1)} \cdots
g_{\pi (k)}$. But this last expression is exactly $\phi ([c])(y)$,
as desired.

For the converse assume $L^\beta G$ is a $G$--graded form of
$\mathbb{C}^c G$.  Then $L^\beta G$ satisfies the graded identities
of $M_n(\mathbb{C})$ and so there is a graded homomorphism from
$U_G$ sending each $x_{ig}$ to $u_g$ in $L^\beta G$.  But under this
map the images of the elements of $M$ are nonzero elements in $L$
and so we get an induced homomorphism from $S^sG$ that takes $S$
into $L$. The result follows.\qed

We can use this result to obtain a parametrization of the
$G$--graded forms over $L$. We begin  with a homomorphism $\varphi$
from $S$ to $L$ and the two-cocycle $\beta$ given by
$\beta(g,h)=\varphi(s(g,h))$.  In the previous proposition we say
that $\beta$ is cohomologous to $c$ over $\mathbb{C}$ and so
$L^\beta G$ is a form of $M_n(\mathbb{C})$. We will now show how to
produce all other forms over $L$.   Let $V$ the free abelian group
(of rank $\ord(G)$) on symbols $r_{g}, g\in G$.  Let $U$ be the
subgroup generated by the elements $r_{g}r_{h}/r_{gh}$. Note that
$V/U$ is a finite group, isomorphic to $G/G^{'}=G_{ab}$ via the map
that sends $r_g$ to $gG'$. For each element $\psi$ in
$\Hom(U,L^\times)$ let $\beta_{\psi}$ be given by
$\beta_{\psi}(g,h)=\beta(g,h)\psi(r_gr_h/r_{gh})$. Note that there
is a canonical homomorphism from $\Hom(V,L^\times)$ to
$\Hom(U,L^\times)$.  We will let $\im((\Hom(V,L^\times))$ denote its
image in $\Hom(U,L^\times)$.

\begin{prop}
\label{parametrization.prop}
The following hold:

(a) For every $\psi\in \Hom(U,L^\times)$ the function $\beta_{\psi}$
is a two-cocycle cohomologous to $\beta$ over $\mathbb{C}$.  In
particular $L^{\beta_{\psi}}G$ is a $G$--graded form of
$M_n(\mathbb{C})$.

(b) If $L^\gamma G$ is a $G$--graded form of $M_n(\mathbb{C})$ there
is a homomorphism $\psi\in \Hom(U,L^\times)$ such that
$\gamma=\beta_{\psi}$. In particular  $L^\gamma G$ is isomorphic to
$L^{\beta_{\psi}}G$ as a $G$--graded algebra.

(c) Two forms $L^{\beta_{\psi}} G$ and $L^{\beta_{\psi '}}G$ are
$G$--graded isomorphic if and only if
$\beta_{\psi}=f(\alpha)\beta_{\psi '}$ for some $\alpha$ in
$\im((\Hom(V,L^\times))$.

(d) The function $\psi\rightarrow \beta_{\psi}$ induces a one-to-one
correspondence between the group
$\Hom(U,L^\times)/\im(\Hom(V,L^\times)$ and the set  of $G$--graded
isomorphism classes of $G$--graded forms of $\mathbb{C}^cG$ over
$L$.
\end{prop}

\proof (a)  It is clear that $\beta_{\psi}$ is a cocycle.  To show
that $\beta$ and $\beta_{\psi}$ are cohomologous over $\mathbb{C}$ it
suffices to show that they determine the same
element in $\Hom(M(G),\mathbb{C}^\times)$, that is that
$\phi([\beta_{\psi}])=\phi([\beta])$, where $\phi$ is the  isomorphism from
$H^{2}(G,\mathbb{C}^\times)$ to $\Hom(M(G),\mathbb{C}^\times)$.  The
argument is quite similar to the proof of the corresponding fact in
Proposition~\ref{forms.prop}, so we will omit it.

(b) If $L^\gamma G$ is a $G$--graded form of $M_n(\mathbb{C})$ then
we have seen that $\gamma$ is cohomologous to $\beta$ over
$\mathbb{C}$ so there is a cochain $\lambda:G\rightarrow
\mathbb{C}^\times$ such that for all $g,h\in G$, we have $\gamma
(g,h)=(\lambda(g)\lambda(h)/\lambda(gh))\beta (g,h)$.  It follows
that for all $g,h\in G$, $\lambda(g)\lambda(h)/\lambda(gh)$ is in
$L$. The cochain $\lambda$ clearly determines a homomorphism from
$V$ to $\mathbb{C}^\times$ (sending $r_g$ to $\lambda (g)$) and this
homomorphism restricts to a homomorphism from $\psi$ from $U$ to
$L^\times$ (given by $\psi
(r_gr_h/r_{gh})=\lambda(g)\lambda(h)/\lambda(gh)$) which satisfies
$\gamma=\beta_{\psi}$.

(c) If $L^{\beta_{\psi}} G$ and $L^{\beta_{\psi '}}G$ are
$G$--graded isomorphic then $\beta_{\psi}$ is cohomologous to
$\beta_{\psi '}$ over $L$.  As in the proof of part (b) this means
there is a homomorphism $\alpha$ from $V$ to $L^\times$ such that
$\beta_{\psi}(g,h)=(\alpha (g)\alpha (h)/\alpha (gh))\beta_{\psi
'}$. It follows that $\psi=\alpha \psi '$, as desired.

(d) This follows from the first three parts.\qed

 We can connect this parametrization with cohomology as follows.  We have a short exact
sequence of groups

$$1\longrightarrow U\longrightarrow V\longrightarrow
G_{ab}\longrightarrow 1$$

Applying the functor $\Hom( - ,L^\times)$ we obtain a long exact
sequence in cohomology

$$\dots\longrightarrow \Hom(V,L^\times)\longrightarrow
\Hom(U,L^\times) \longrightarrow Ext(G_{ab},L^\times)\longrightarrow
1$$ where the surjectivity of the last map follows from the fact
that $Ext(V,L^\times)=0$ because $V$ is a free abelian group.  We
also have the universal coefficients sequence:

$$1\rightarrow Ext^{1}(G_{ab},L^\times)
\rightarrow H^{2}(G,L^\times)\rightarrow
\Hom(M(G),L^\times)
\rightarrow 1$$

Combining these equations shows that
$\Hom(U,L^\times)/\im(\Hom(V,L^\times)$ is isomorphic to
$Ext^{1}(G_{ab},L^\times)$, the kernel of the map from
$H^{2}(G,L^\times)$ to $\Hom(M(G),L^\times)$.

The last result in this section is the determination of the rank of
$Y_f$, the free part of the subgroup of $k^\times$ generated by
the values of $s$.  The computation involves the groups $U$ and $V$.
We first remark that because $V/U$ is finite, the rank of $U$ is equal
to the rank of $V$, which is clearly $n$, the order of $G$.

\begin{prop}
\label{rank.prop} The group $Y_f$ is isomorphic to $U$. In
particular the rank of $Y_f$ is $n$, the order of $G$.
\end{prop}

\proof Recall that $Y$ denotes the subgroup of $k^\times$ generated
by the values of the cocycle $s$.  Let $H$ denote the subgroup of
$k^\times$ generated by the values of the cocycle $c$.  Then by
Proposition~\ref{fieldF.prop} $Y\cap H= Y\cap
\mathbb{C}^\times=\mu$.  Let $\tilde{U}$ be the subgroup of
$k^\times$ generated by set $\{t_gt_h/t_{gh} \ | \ g,h\in G\}$.
Clearly the group $\tilde{U}$ is isomorphic to $U$. Now we have
$Y_f\cong Y/{Y_t} \cong Y/{Y\cap H}\cong {YH}/H$.  But
$YH=\tilde{U}H$.  Hence ${YH}/H={\tilde{U}H}/H\cong
\tilde{U}/{\tilde{U}\cap H}\cong \tilde{U}$ because $\tilde{U}\cap
H=1$. Hence $Y_f\cong\tilde{U}\cong U$, as desired.\qed

\begin{cor}
\label{Srank.cor}
The field of fractions of $S$ has transcendence degree $n$ over
$\mathbb{Q}(\mu)$.
\end{cor}

\proof This follows from the description of $S$ given in
Proposition~\ref{fieldF.prop}.\qed

\bigskip

\section{The universal $G$--graded central simple algebra.}

We have seen that each fine grading by a group $G$ on
$M_n(\mathbb{C})$ gives rise to a nondegenerate two-cocycle $c$ on
the group $G$ (which must be of central type) such that
$M_n(\mathbb{C})$ is graded isomorphic to $\mathbb{C}^cG$.  We have
also seen that the universal $G$--graded algebra $U_G$ is a prime
ring with ring of central quotients $Q(U_G)$ graded isomorphic to
the $F$--central simple algebra $F^s G$ where $s$ is the generic
cocycle obtained from $c$ and $F=\mathbb{Q}(\mu)(y_1,\dots,
y_m)({t_{ig} \over t_g} \ | \  i\geq 1, g\in G)$ (see
Proposition~\ref{fieldF.prop}). In this section we obtain
information about the index of $Q(U_G)$ and about the dependence of
$Q(U_G)$ on the given cocycle $c$.

\begin{prop} \label{indexbound.prop} The index of $Q(U_G)$ is equal to
the maximum of the indices of the $G$--graded forms $L^\beta G$
where $L$ varies over all subfields of $\mathbb{C}$ that contain
$\mathbb{Q}(\mu)$ and $\beta$ varies over all nondegenerate cocycles
over $L$ cohomologous to $c$ over $\mathbb{C}$.
\end{prop}

\proof Let $M$ denote the maximum of the indices of the forms. We
first claim that $\ind(Q(U_G))\leq M$: In fact we can specialize the
universal algebra $U_G$ by sending $x_g$ to $r_g u_g$ where the
elements $\{r_g \mid g\in G\}$ are algebraically independent complex
numbers.  The resulting algebra $\overline{U}$ is isomorphic to
$U_G$ and tensoring with the field of fractions of its center gives
a central simple algebra that is isomorphic to $Q(U_G)$.  Clearly
the index of this form equals the index of  $Q(U_G)$ and so we get
the inequality. We proceed to show that the index of $Q(U_G)=F^sG$
is greater than or equal to $M$. Recall
(Proposition~\ref{fieldF.prop}) that the center $S$ of $S^{s}G$ is a
ring of Laurent polynomials over the field $\mathbb{Q}(\mu)$. If $M$
is any maximal ideal of $S$ then $S_M$ is a regular local ring.  By
Proposition~\ref{forms.prop} it suffices to show that if $M$ is any
maximal ideal of $S$ then the index of the algebra $S^sG\otimes_S F=
F^{s}G$ is greater than or equal to the index of the residue algebra
$S^sG/MS^sG$.  This is a well known fact.  We outline a proof: Let
$r_1,r_2,\dots,r_m$ be a system of parameters of the regular local
ring $R=S_M$ and let $F^{s}G=M_t(D)$ where $D$ is an $F$--central
division algebra (so $\ind(F^{s}G) = \ord(G)^{1/2}/t$). The
localization $R_{(r_1)}$ is a discrete valuation ring with residue
field the field of fractions of the regular local ring
$R_1=R/(r_1)$. The ring $R_{(r_1)}^sG$ is an Azumaya algebra over
the discrete valuation ring $R_{(r_1)}$ and hence is isomorphic to
$M_t(A)$ where $A$ is an Azumaya algebra such that $A\otimes F=D$
(see Reiner \cite[Theorem 21.6]{R}). It follows that $R_1^sG$ is an
Azumaya algebra and if we let $Q(R_1)$ denote the field of fractions
of $R_1$, then $R_1^s G\otimes Q(R_1)$ is isomorphic to $M_t(D_1)$
where $D_1$ a central simple $Q(R_1)$--algebra, the residue algebra
of $D$. The image of $M$ in $R_1$ is generated by the images of
$r_2,r_3,\dots, r_m$  and these $m-1$ elements form a system of
parameters for $R_1$. We can therefore repeat the process and
eventually obtain that $S_M^sG/MS_M^sG$ is isomorphic to an algebra
of the form $M_t(E)$ for some central simple $S/M$--algebra $E$.
Hence the index of $S_M^sG/MS_M^sG$ is at most the index of $F^sG$.
\qed

The question for which groups $G$ of central type there is a cocycle
$c$ such that $Q(U_{G,c})$ is a division algebra was answered in
Aljadeff, et al \cite{AHN} and Natapov \cite{N}: We consider the
following list of $p$-groups, called $\Lambda_p$:

\begin{enumerate}
   \item $G$
   is abelian of symmetric type, that is $G \cong \prod(C_{p^{n_i}} \times C_{p^{n_i}})$,
   \item $G \cong G_1 \times G_2$ where \\ $G_1 = C_{p^n} \rtimes C_{p^n} =
   \langle\pi, \sigma \mid \sigma^{p^n}=\pi^{p^n}=1 \; {\rm and} \;
   \sigma\pi\sigma^{-1} = \pi^{p^s+1}\rangle$ where $1 \leq s < n$ and
   $1 \neq s$ if $p=2$,
   and \\ $G_2$ is an abelian group of symmetric type of exponent
   $\leq p^s$,
   \item $G \cong G_1 \times G_2$ where \\
   $G_1 = C_{2^{n+1}} \rtimes (C_{2^n} \times Z_2) =
               \bigg \langle \pi,\sigma,\tau \, \bigg|  \,
               \begin{array}{l}
                 \pi^{2^{n+1}}=\sigma^{2^n}=\tau^2=1, \sigma\tau = \tau\sigma,\\
                 \sigma\pi\sigma^{-1}=\pi^3,
                 \tau\pi\tau^{-1}=\pi^{-1}
               \end{array}\bigg\rangle$ and \\
   $G_2$ is an abelian group of symmetric type of exponent $\leq 2$.
\end{enumerate}

We let $\Lambda$ be the collection of nilpotent groups such that for
any prime $p$, the Sylow-$p$ subgroup is on the $\Lambda_p$.

\begin{prop}
\label{list.prop}
Let $G$ be a group of central type of order $n^2$.  The following
conditions are equivalent:

(a) There is a nondegenerate cocycle $c$ on $G$ for which the
universal $G$--graded algebra $U_{G,c}$ is a domain.

(b) There is a nondegenerate cocycle $c$ on $G$ for which the
universal central simple algebra $Q(U_{G,c})$ is a division algebra.

(c) There is a nondegenerate cocycle $c$ on $G$ such that the
resulting $G$--graded algebra $M_n(\mathbb{C})$ has a $G$--graded
form that is a division algebra.

(d) The group $G$ is on the list $\Lambda$.
\end{prop}

\proof The equivalence of (a) and (b) follows from the fact that
$F^sG$ is isomorphic to $Q(U_{G,c})$. The equivalence of (b) and (c)
follows from the previous proposition. Finally the equivalence of
(c) and (d) follows from Aljadeff et al \cite[Corollary 3]{AHN} and
Natapov \cite[Theorem 3]{N}. \qed

Our main result in this section is that if $G$ is a group on the
list $\Lambda$ then the universal central simple $G$--graded algebra
$Q(U_{G,c})$ is independent (up to a non-graded isomorphism) of the
cocycle $c$. In fact we will show that for groups $G$ on the list
the automorphism group of $G$ acts transitively on the set of
classes of nondegenerate cocycles.

To begin let $G$ be any central type group (not necessarily on
$\Lambda$) and let $c$ be a nondegenerate two-cocycle on $G$ with
coefficients in $\mathbb{C}^\times$. Let $\varphi$ be an
automorphism of $G$ and let $\varphi(c)$ be the two-cocycle defined
by $\varphi(c)(\sigma,\tau) = c (\varphi^{-1}(\sigma),
\varphi^{-1}(\tau))$. It is clear that $\varphi(c)$ is also
nondegenerate.
\begin{thm}
\label{transitive.thm} If $G$ is a group on the list $\Lambda$ and
$c$ and $c'$ are nondegenerate cocycles on $G$ with values in
$\mathbb{C}^\times$, then there is an automorphism $\varphi$ of $G$
such that $\varphi(c)$ is cohomologous over $\cc$ to $c'$.
\end{thm}

\proof Let $G$ be a group on the list $\Lambda$. Because $G$ is
necessarily nilpotent we can assume that $G$ lies on $\Lambda_p$ for
some prime $p$.  In the course of the proof whenever we refer to
basis elements $\{u_{g}\}_{g\in G}$ in the twisted group algebra
$\cc^c G$, we assume they satisfy $u_g u_h = c(g,h) u_{gh}$.

The strategy is as follows: For each group $G \in \Lambda_{p}$ we
fix a set of generators $\Phi$. Then we consider the family of sets
of generators $\Phi^{'}=\varphi(\Phi)$ that are obtained from $\Phi$
via an automorphism $\varphi$ of $G$. We say that $\Phi$ and
$\Phi^{'}$ are equivalent. Next, we exhibit a certain nondegenerate
cohomology class $\alpha \in H^2(G, \cc^\times)$ by setting a set of
equations (denoted by $E_{G}$) satisfied by elements
$\{u_{g}\}_{g\in \Phi}$ in $\mathbb{C}^{c_{0}}G$, where $c_0$ is a
two-cocycle representing $\alpha$. We say that the cohomology class
$\alpha$ (or by abuse of language, the two-cocycle $c_0$) is of
``standard form" with respect to $\Phi$. It is indeed abuse of
language since in general the equations do not determine the
two-cocycle $c_{0}$ uniquely but only its cohomology class. Finally
we show that for any nondegenerate two-cocycle $c$ with values in
$\mathbb{C}^\times$ there is a set of generators
$\Phi^{'}=\varphi(\Phi)$ with respect to which $c$ is of standard
form. The desired automorphism of the group $G$ will be determined
by compositions of ``elementary" automorphisms, that is
automorphisms that are defined by replacing some elements of the
generating set. We start with $\Phi^{(0)}=\Phi$ and denote the
updated generating sets by $\Phi^{(r)}=\{g_{1}^{(r)},\ldots
,g_{n}^{(r)}\}, r=1,2,3,\dots$. Elements in the generating set that
are not mentioned remain unchanged, that is
$g_{i}^{(r)}=g_{i}^{(r-1)}$. In all steps it will be an easy check
that the map defined is indeed an automorphism of $G$.

{\bf (I)} We start with abelian groups on $\Lambda_{p}$. Let
$G=C_{p^{r_{1}}}\times C_{p^{r_{1}}}\times C_{p^{r_{2}}}\times
C_{p^{r_{2}}}\times \dots \times C_{p^{r_{m}}} \times C_{p^{r_{m}}}$
and let $\Phi = \{ \gamma_1, \gamma_2, \dots, \gamma_{2m} \}$ be an
ordered set of generators.

Fixing a primitive $p^n$-th root of unity $\varepsilon$, where $n
\geq \max(r_k)$, $1\leq k \leq m$, we say a two-cocycle $c_{0}$ is
of standard form with respect to $\Phi$ if there are elements
$\{u_{g}\}_{g\in \Phi}$ in $\mathbb{C}^{c_0}G$ that satisfy

\begin{equation}\label{abelian_normal_form}
\left\{
    \begin{array}{ll}
(u_{\gamma_{2k-1}},u_{\gamma_{2k}})=\varepsilon^{p^{n-r_k}} \ \
\hbox{for all} \ 1\leq k\leq m, \\
\hbox{all other commutators of} \ {u_{\gamma_i}}'s \ \hbox{are
trivial.}
    \end{array}
  \right.
\end{equation}

Note that $[c_0]$ is determined by (\ref{abelian_normal_form}).
Given any nondegenerate two-cocycle $c$ it is known that there is a
set of generators $\Phi^{'}=\varphi(\Phi)$ with respect to which $c$
is of standard form (see, for example, Aljadeff and Sonn
\cite[Theorem 1.1]{AS}). Hence we are done in this case.

\medskip {\bf (II)} Next we consider the group $G = C_{p^n}\rtimes C_{p^n} \times C_{p^{r_1}} \times
C_{p^{r_1}} \times \dots \times C_{p^{r_m}} \times C_{p^{r_m}}$ with
a set of generators $\Phi = \{ \pi,\sigma, \gamma_1, \dots,
\gamma_{2m}\}$ where $\pi, \sigma$ satisfy
$\sigma(\pi)=\pi^{p^{s}+1}$ for some $s=1,\dots,n-1$ if $p$ is odd
and for some $s=2,\dots,n-1$ if $p=2$, and  $r_k \leq s$ for all
$1\leq k \leq m$. We write $G = G_1 \times G_2$ where $G_1 = C_{p^n}
\rtimes C_{p^n}$.

Fix a primitive $p^{n}$-th root of unity $\varepsilon$. There exist
a two-cocycle  $c_1 \in Z^2(G_1, \cc^\times)$ and elements
$\{{u_g}\}_{g \in \{\pi,\sigma\}}$ in $\cc^{c_1} G_1$ that satisfy
\begin{equation}\label{p odd_normal_form}
u_\pi^{p^n}=u_\sigma^{p^n}=1, u_\sigma u_\pi u_\sigma^{-1} =
\varepsilon u_\pi^{p^{s+1}}.
\end{equation}
Moreover the class $[c_1] \in H^2(G_1, \cc^\times)$ is uniquely
determined by (\ref{p odd_normal_form}) (see, for example, the proof of
Karpilovsky \cite[Theorem 10.1.25]{Ka}).

Now, by Karpilovsky (\cite[Proposition 10.6.1]{Ka}) we have $ M(G)
\cong M(G_1) \times M(G_2) \times ( G_1^{ab} \otimes G_2^{ab})$
(here $G^{ab}$ denotes the abelianization of $G$) and hence there
exist a two-cocycle $c_0 \in Z^2(G, \cc^\times)$ and elements
$\{{u_g}\}_{g \in \Phi}$ in $\cc^{c_0} G$ that satisfy equations
(\ref{p odd_normal_form}), (\ref{abelian_normal_form}) and
$u_{\pi}u_{\gamma_{i}}=u_{\gamma_{i}}u_{\pi}$,
$u_{\sigma}u_{\gamma_{i}}=u_{\gamma_{i}}u_{\sigma}$. We denote this
set of equations by $E_{G}$. The class $[c_0] \in H^2(G,
\cc^\times)$ is uniquely determined by $E_{G}$ and is easily seen to
be nondegenerate. We say that a class $[c_0]$ is of standard form
with respect to $\Phi$ if there are elements $\{{u_g}\}_{g \in
\Phi}$ in $\cc^{c_0} G$ that satisfy $E_{G}$.

Now, let $c\in Z^2(G, \cc^\times)$ be any nondegenerate two-cocycle,
and let $u_g$ be representatives of elements of $G$ in $\cc^{c} G$.
As explained above we set $\Phi^{(0)}=\{\pi^{(0)}, \sigma^{(0)},
\gamma^{(0)}_1, \dots, \gamma^{(0)}_{2m}\}=\Phi=\{\pi, \sigma,
\gamma_1, \dots, \gamma_{2m}\}$. We are to exhibit a sequence of
automorphisms such that their composition applied to $\Phi^{(0)}$
yields a generating set with respect to which the cocycle $c$ is of
standard form. First, note that we may assume (by passing to an
equivalent cocycle, if necessary) that $u_{\pi^{(0)}}^{p^n} =
u_{\sigma^{(0)}}^{p^n} = 1$. Let $\alpha \in \cc^\times$ be
determined by $u_{\sigma^{(0)}} u_{\pi^{(0)}} u_{\sigma^{(0)}}^{-1}
= \alpha u_{\pi^{(0)}}^{p^s+1}$. It follows that $u_{\sigma^{(0)}}
u_{\pi^{(0)}}^{p^{n-1}} u_{\sigma^{(0)}}^{-1} = \alpha^{p^{n-1}}
u_{\pi^{(0)}}^{p^{n-1}}$ and hence, $\alpha^{p^{n-1}}$ is a $p$-th
root of unity. Next, observe that since $c$ is nondegenerate,
$u_{\sigma^{(0)}}$ cannot commute with $u_{\pi^{(0)}}^{p^{n-1}}$
(for otherwise, $u_{\pi^{(0)}}^{p^{n-1}}$ is contained in the center
of the algebra) and hence $\alpha$ is a primitive $p^n$-th root of
unity. Now, replacing $\pi^{(0)}$ by a suitable prime to $p$ power
of $\pi^{(0)}$ (that is $\pi^{(1)}=(\pi^{(0)})^{l}$, $l$ prime to
$p$), we may assume that $\alpha = \varepsilon$. Leaving
$\sigma^{(0)}$ and all other generators unchanged, that is
$\sigma^{(1)}=\sigma^{(0)}, \gamma^{(1)}_1=\gamma^{(0)}_1, \dots,
\gamma^{(1)}_{2m}=\gamma^{(0)}_{2m}$, we obtain a generating set
$\Phi^{(1)} = \{\pi^{(1)},\sigma^{(1)}, \gamma^{(1)}_1, \dots,
\gamma^{(1)}_{2m} \}$ of $G$, equivalent to $\Phi$, such that the
elements $u_{\pi^{(1)}},u_{\sigma^{(1)}}$ in $\cc^{c} G_1$ satisfy
(\ref{p odd_normal_form}).

Assume $(u_{\gamma^{(1)}_i}, u_{\pi^{(1)}}) = \xi_i \neq 1$ (clearly
$\xi_i$ is a root of unity of order $p^r (\leq \ord(\gamma^{(1)}_i)
\leq p^s))$. Now, by the nondegeneracy of $c$, the root of unity
$\zeta$ determined by the equation $u_{\sigma^{(1)}}^{p^{n-s}}
u_{\pi^{(1)}} u_{\sigma^{(1)}}^{-p^{n-s}} = \zeta u_{\pi^{(1)}}$ is
of order $p^s$ and hence $\xi_i = \zeta^l$ for some $l$. Observe
that $\ord((({\sigma^{(1)}})^{p^{n-s}})^l) \leq p^r$ and hence we
may put $\gamma^{(2)}_i= \gamma^{(1)}_i ({\sigma^{(1)}})^{-l
p^{n-s}}$ and get $(u_{\gamma^{(2)}_i}, u_{\pi^{(1)}}) = 1$.
Performing this process for all $1\leq i \leq 2m$ (and leaving
$\pi^{(1)}$ and $\sigma^{(1)}$ unchanged) we get a set of generators
$\Phi^{(2)}$ such that $(u_{\gamma^{(2)}_i}, u_{\pi^{(2)}}) = 1$ for
all $i$.

Next, if $(u_{\gamma^{(2)}_i}, u_{\sigma^{(2)}}) \neq 1$ we set
$\gamma^{(3)}_i=\gamma^{(2)}_i ({\pi^{(2)}})^{-t p^{n-s}}$ for an
appropriate positive integer $t$ and get $(u_{\gamma^{(3)}_i},
u_{\sigma^{(2)}}) = 1$. Performing this process for all $1\leq i
\leq 2m$ we get a set of generators $\Phi^{(3)}$ such that
$(u_{\gamma^{(3)}_i}, u_{\sigma^{(3)}}) = 1$ for all $i$.

Finally, we may proceed as in Case I and obtain a generating set
$\Phi^{(4)}$, equivalent to $\Phi$, with respect to which $c$ is of
standard form.

\medskip {\bf (III)} Next we consider the group
$G = C_{2^{n+1}}\rtimes (C_{2^n}\times C_2)\times C_2 \times C_2
\ldots \times C_2 \times C_2$, $n\geq 2$, with a set of generators
$\Phi=\{\pi,\sigma,\tau, \gamma_1, \dots, \gamma_{2m}\}$, where
$\sigma\pi\sigma^{-1}=\pi^3$ and $\tau\pi\tau^{-1}=\pi^{-1}$. We
write $G = G_1 \times G_2$ where $G_1 = C_{2^{n+1}}\rtimes
(C_{2^n}\times C_2)$ and $G_2$ is elementary abelian of rank $2m$.

We first exhibit a construction of a field $K$ and a cocycle $\beta
\in Z^2(G_1, K^\times)$ such that the algebra $D \cong K^{\beta}G_1$
is $K$--central simple (in fact a $K$--central division algebra),
which is analogous to that in Natapov \cite{N}. Let $K = \qq(s,t)$
be the subfield of $\cc$ generated by algebraically independent
elements $s$ and $t$. Let $L = K(v_\pi)/K$ be a cyclotomic extension
where $v_\pi^{2^{n+1}} = -1$. The Galois action of $Gal(L/K) \cong
Z_{2^n} \times Z_2 = \langle \sigma, \tau \rangle$ on $L$ is given
by
\[
 \sigma(v_\pi) =  v_\pi^3 \;\;\; {\rm and} \;\;\; \tau(v_\pi) =
 v_\pi^{-1}.
\]
The algebra $D$ is the crossed product $D=(L/K, \langle \sigma, \tau
\rangle)$ determined by the following relations
\[
 v_\sigma^{2^n} = s, \; v_\tau^2 = t, \;  (v_\sigma, v_\tau) = 1,
\]
where $v_\sigma$ and $v_\tau$ represent $\sigma$ and $\tau$ in $D$.
It is easy to see that $D$ is isomorphic to a twisted group algebra
of the form $K^{\beta} G_1$.

Fix a primitive $2^{n+1}$-th root of unity $\varepsilon$, and let
$\xi$ be a square root of $\varepsilon$. Let $\widetilde{s}$ and
$\widetilde{t}$ be elements of $\cc$ such that $\widetilde{s}^{2^n}
= s^{-1}$ and $\widetilde{t}^2 = t^{-1}$. Then the elements
$u_{\pi^i \sigma^j \tau^k} = (\xi v_\pi)^i (\widetilde{s}
v_\sigma)^j (\widetilde{t} v_\tau)^k$ in $\cc^{\beta} G_1 \cong
\cc\otimes_{K}D$ satisfy

\begin{equation}\label{p 2_normal_form}
u_\pi^{2^{n+1}}=u_\sigma^{2^n}=u_\tau^2=1, \ u_\sigma u_\pi =
\varepsilon u_\pi^3 u_\sigma, \ u_\tau u_\pi = \varepsilon^{-1}
u_\pi^{-1}u_\tau, \ u_\tau u_\sigma =u_\sigma u_\tau.
\end{equation}

Clearly, there exists a two-cocycle $c_1 \in Z^2(G_1, \cc^\times)$,
cohomologous to $\beta$, such that $u_{g_1}u_{g_2} =
c_1(g_1,g_2)u_{g_1g_2}$ for all $g_1, g_2 \in G_1$  and moreover the
class $[c_1]$ is determined by (\ref{p 2_normal_form}).

Now, by Karpilovsky \cite[Proposition 10.6.1]{Ka} there exist a
two-cocycle $c_0 \in Z^2(G, \cc^\times)$ and elements $\{{u_g}\}_{g
\in \Phi}$ in $\cc^{c_0} G$ that satisfy (\ref{p 2_normal_form}),
(\ref{abelian_normal_form}) and
$u_{g}u_{\gamma_{i}}=u_{\gamma_{i}}u_{g}$ for $g\in \{\pi, \sigma,
\tau\}$. We denote this set of equations by $E_{G}$. As in the
previous cases, also here, $[c_0]$ is determined by $E_G$, and we
say that the cocycle $c_{0}$ is of standard form with respect to
$\Phi$.

Now, let $c\in Z^2(G, \cc^\times)$ be a nondegenerate two-cocycle,
and consider the algebra $\cc^{c} G = \oplus_{g \in G} \cc u_g$. Let
$\Phi^{(0)}= \{\pi^{(0)},\sigma^{(0)},\tau^{(0)}, \gamma^{(0)}_1,
\dots, \gamma^{(0)}_{2m}\} = \Phi=\{\pi,\sigma,\tau, \gamma_1,
\dots, \gamma_{2m}\}$. We may assume (by passing to an equivalent
cocycle, if necessary) that $u_{\pi^{(0)}}^{2^{n+1}} =
u_{\sigma^{(0)}}^{2^n} = u_{\tau^{(0)}}^{2} = 1$. Let $\alpha \in
\cc^\times$ be determined by the equation $u_{\sigma^{(0)}}
u_{\pi^{(0)}} u_{\sigma^{(0)}}^{-1} = \alpha u_{\pi^{(0)}}^3$. A
straightforward calculation shows that for any $k$ in $\nn$,
\begin{equation}\label{sigmaction}
u_{\sigma^{(0)}}^k u_{\pi^{(0)}} u_{\sigma^{(0)}}^{-k} =
\alpha^{\frac{3^k-1}{2}}u_{\pi^{(0)}}^{3^k}, \  2^k \mid
\frac{3^{2^k-1}}{2} \ \hbox{and} \  2^{k+1} \nmid
\frac{3^{2^k-1}}{2}.
\end{equation}

We \underline{claim} that $\alpha$ is a primitive $2^{n+1}$-th root
of unity. To see this note that $\sigma^{(0)}$ and
$(\pi^{(0)})^{2^n}$ commute and therefore, $u_{\sigma^{(0)}}
u_{\pi^{(0)}}^{2^n} u_{\sigma^{(0)}}^{-1} = \pm
u_{\pi^{(0)}}^{2^n}$. We need to show that $u_{\sigma^{(0)}}
u_{\pi^{(0)}}^{2^n} u_{\sigma^{(0)}}^{-1} = -u_{\pi^{(0)}}^{2^n}$.
If not we have
$$
u_{\pi^{(0)}}^{2^n} = (u_{\sigma^{(0)}} u_{\pi^{(0)}}
u_{\sigma^{(0)}}^{-1})^{2^n} = (\alpha u_{\pi^{(0)}}^3)^{2^n} =
\alpha^{2^n} u_{\pi^{(0)}}^{2^n},
$$
that is $\alpha^{2^n} = 1$. Then, by (\ref{sigmaction}) we have
$$
u_{\sigma^{(0)}}^{2^{n-1}} u_{\pi^{(0)}} u_{\sigma^{(0)}}^{-2^{n-1}}
= \alpha^{\frac{3^{2^{n-1}}-1}{2}}u_{\pi^{(0)}}^{3^{2^{n-1}}} =
u_{\pi^{(0)}}
$$
and hence $u_{\sigma^{(0)}}^{2^{n-1}}$ is in the center of $\cc^{c}
G$. This contradicts the nondegeneracy of the cocycle $c$ and the
claim is proved. It follows that if we replace $\pi^{(0)}$ by a
suitable odd power of $\pi^{(0)}$ (that is
$\pi^{(1)}={\pi^{(0)}}^{l}$, $l$ is odd) and leaving the other
generators unchanged we obtain a generating set $\Phi^{(1)} =
\{\pi^{(1)},\sigma^{(1)}, \tau^{(1)} \gamma^{(1)}_1, \dots,
\gamma^{(1)}_{2m} \}$ of $G$, equivalent to $\Phi$, such that
$u_{\sigma^{(1)}} u_{\pi^{(1)}} u_{\sigma^{(1)}}^{-1} = \varepsilon
u_{\pi^{(1)}}^3$.

Next, since $\sigma^{(1)}$ and $\tau^{(1)}$ commute and $\tau^{(1)}$
is of order 2 we have that $(u_{\sigma^{(1)}}, u_{\tau^{(1)}}) = \pm
1$. If the commutator is $-1$, we put $\tau^{(2)}=\tau^{(1)}
{\pi^{(1)}}^{2^n}$ and $\sigma^{(2)}=\sigma^{(1)}$. Otherwise we
leave both unchanged and, in either case, get an equivalent
generating set $\Phi^{(2)}$ with $(u_{\sigma^{(2)}}, u_{\tau^{(2)}})
= 1$.

Let $\beta \in \cc^\times$ be determined by the equation
$u_{\tau^{(2)}} u_{\pi^{(2)}} u_{\tau^{(2)}}^{-1} = \beta
u_{\pi^{(2)}}^{-1}$. Then
$$
u_{\sigma^{(2)}} u_{\tau^{(2)}} u_{\pi^{(2)}} u_{\tau^{(2)}}^{-1}
u_{\sigma^{(2)}}^{-1}= \beta \varepsilon^{-1} u_{\pi^{(2)}}^{-3}, \
\hbox{and} \ u_{\tau^{(2)}} u_{\sigma^{(2)}} u_{\pi^{(2)}}
u_{\sigma^{(2)}}^{-1} u_{\tau^{(2)}}^{-1} = \varepsilon \beta^3
u_{\pi^{(2)}}^{-3}
$$
and by equality of the left hand sides we have that $\beta = \pm
\varepsilon^{-1}$. Now put $\tau^{(3)}=\tau^{(2)}
(\sigma^{(2)})^{2^{n-1}}$ if $\beta = - \varepsilon^{-1}$ and
$\tau^{(3)}=\tau^{(2)}$ otherwise. Leaving the other generators
unchanged we obtain an equivalent set of generators $\Phi^{(3)}$
such that the elements $u_{\tau^{(3)}}, u_{\pi^{(3)}} \in \cc^{c} G$
satisfy $u_{\tau^{(3)}} u_{\pi^{(3)}} u_{\tau^{(3)}}^{-1} = - \beta
u_{\pi^{(3)}}^{-1} = \varepsilon^{-1} u_{\pi^{(3)}}^{-1}$. Thus we
have shown there exist $\{{u_{g}}\}_{g \in
\{\pi^{(3)},\sigma^{(3)},\tau^{^{(3)}}\}}$ in $\cc^{c} G_1$ that
satisfy (\ref{p 2_normal_form}).

Next, we may proceed as in Aljadeff, et al (\cite[Proposition 13]{AHN}) to get a
generating set $\Phi^{(4)}$ of $G$, equivalent to $\Phi$, such that
(a) (\ref{p 2_normal_form}) is satisfied (b)
$(u_{\gamma^{(4)}_{2k-1}},u_{\gamma^{(4)}_{2k}})= \pm 1$  for all
$1\leq k\leq m$ and (c) all other commutators of the
$u_{\gamma^{(4)}_i}$'s are trivial.

The next step is to add the condition $u_{g} u_{\gamma_i} =
u_{\gamma_i} u_{g}$ for $g\in \{\pi, \sigma, \tau \}$:

For each $\gamma^{(4)}_i \in \{\gamma^{(4)}_1, \dots,
\gamma^{(4)}_{2m}\}$, if $(u_{\gamma^{(4)}_i}, u_{\pi^{(4)}}) = - 1$
then put $\gamma^{(5)}_i = \gamma^{(4)}_i (\sigma^{(4)})^{2^{n-1}}$,
otherwise leave $\gamma^{(4)}_i$ unchanged. Leaving $\pi^{(4)},
\sigma^{(4)}, \tau^{(4)}$ unchanged we obtain a generating set
$\Phi^{(5)}$ such that $(u_{\gamma^{(5)}_i}, u_{\pi^{(5)}}) =1$ for
all $i$. Next, if $(u_{\gamma^{(5)}_i}, u_{\sigma^{(5)}}) = - 1$ put
$\gamma^{(6)}_i = \gamma^{(5)}_i (\pi^{(5)})^{2^n}$, otherwise leave
$\gamma^{(5)}_i$ unchanged. We obtain a generating set $\Phi^{(6)}$
such that $(u_{\gamma^{(6)}_i}, u_{\sigma^{(6)}}) = 1$ for all $i$.

At this point we have for each $1\leq k \leq  m$,
$(u_{\gamma^{(6)}_{2k-1}},u_{\gamma^{(6)}_{2k}})= \pm 1$. In fact we
\underline{claim} that
$(u_{\gamma^{(6)}_{2k-1}},u_{\gamma^{(6)}_{2k}})= - 1$. If not, by
the preceding steps we have that the element
$u_{\gamma^{(6)}_{2k-1}}$ centralizes all the $u_{\gamma^{(6)}_i}$'s
as well as $u_{\pi^{(6)}}$ and $u_{\sigma^{(6)}}$. It follows that
$u_{\gamma^{(6)}_{2k-1}}$ does not centralize $u_{\tau^{(6)}}$, for
otherwise we get a contradiction to the nondegeneracy of the
two-cocycle $c$. A similar argument shows that
$u_{\gamma^{(6)}_{2k}}$ does not centralize $u_{\tau^{(6)}}$ , but
then the product $u_{\gamma^{(6)}_{2k-1}}u_{\gamma^{(6)}_{2k}}$ is
central, a contradiction. This proves the claim. We refer to the
pair of elements ${\gamma_{2k-1}}, {\gamma_{2k}}$ as {\it partners}.
We proceed now as follows: If $(u_{\gamma^{(6)}_{1}},
u_{\tau^{(6)}}) = - 1$, put $\tau^{(7)} = \tau^{(6)} \gamma^{(6)}_2$
(i.e. multiply $\tau^{(6)}$ by the partner of $\gamma^{(6)}_1$), and
leave $u_{\gamma^{(6)}_1}$ and all the other generators unchanged.
We continue in a similar way with $\gamma^{(6)}_2, \gamma^{(6)}_3,
\dots, \gamma^{(6)}_{2m}$.

We now have a generating set $\Phi^{(r)}$ (some $r$), equivalent to
$\Phi$, such that $\{{u_{g}}\}_{g \in \Phi^{(r)}}$ in $\cc^c G$
satisfy (a) equations (\ref{p 2_normal_form}) (b)
$(u_{\gamma^{(r)}_{2k-1}},u_{\gamma^{(r)}_{2k}})= - 1$  for all
$1\leq k \leq  m$ (c) all other commutators of the
$u_{\gamma^{(r)}_i}$'s are trivial and (d) $u_{\gamma^{(r)}_i}$
centralizes the subalgebra $\cc(u_{\pi^{(r)}}, u_{\sigma^{(r)}},
u_{\tau^{(r)}})$ for all $i$. It follows that the cocycle $c$ is of
standard form with respect to $\Phi^{(r)}$.

\medskip {\bf (IV)} In the last step we consider the group
$G = C_4 \rtimes (C_2 \times C_2) \times C_2 \times C_2 \ldots
\times C_2 \times C_2$ with a set of generators $\Phi =
\{\pi,\sigma,\tau, \gamma_1, \dots, \gamma_{2m}\}$, where
$\sigma\pi\sigma^{-1}=\pi^3$ and $\tau\pi=\pi\tau$. We write $G =
G_1 \times G_2$ where $G_1 = C_4 \rtimes (C_2 \times C_2)$ and $G_2$
is elementary abelian of rank $2m$.

As in the previous case, a construction of a $K$--central division
algebra of the form $K^{\beta} G_1$, analogous to that in Natapov
\cite{N}, yields a $4\times4$ matrix algebra isomorphic to
$\cc^{c_1} G_1$ with a basis $\{u_g\}_{g \in G_1}$ such that
$u_{\pi^i \sigma^j \tau^k} = u_\pi^i u_\sigma^j u_\tau^k$ and
\begin{equation}\label{422_normal_form}
u_\pi^4=u_\sigma^2=u_\tau^2=1, \ u_\sigma u_\pi = \varepsilon
u_\pi^3 u_\sigma, \ u_\tau u_\pi = - u_\pi u_\tau, \ u_\tau u_\sigma
=u_\sigma u_\tau,
\end{equation}
where $\varepsilon$ is a 4-th root of unity. Clearly, the class
$[c_1]$ is determined by (\ref{422_normal_form}).

Now, by Karpilovsky \cite[Proposition 10.6.1]{Ka} there exist a
two-cocycle $c_0 \in Z^2(G, \cc^\times)$ and elements $\{{u_g}\}_{g
\in \Phi}$ in $\cc^{c_0} G$ that satisfy (\ref{422_normal_form}),
(\ref{abelian_normal_form}), and
$u_{g}u_{\gamma_{i}}=u_{\gamma_{i}}u_{g}$ for $g\in \{\pi, \sigma,
\tau\}$. We denote these equations by $E_{G}$. Note that $[c_0]$ is
determined by $E_G$. As usual we will say the cocycle $c_{0}$ is of
standard form with respect to $\Phi$.

Now, let $c\in Z^2(G, \cc^\times)$ be a nondegenerate two-cocycle
and write $\cc^{c} G = \oplus_{g \in G} \cc u_g$. As in previous
steps put $\Phi^{(0)}=\{\pi^{(0)},\sigma^{(0)},\tau^{(0)},
\gamma^{(0)}_1, \dots, \gamma^{(0)}_{2m}\}=\Phi = \{\pi,\sigma,\tau,
\gamma_1, \dots, \gamma_{2m}\}$. We may assume (by passing to an
equivalent cocycle, if necessary) that
$u_{\pi^{(0)}}^4=u_{\sigma^{(0)}}^2=u_{\tau^{(0)}}^2=1$.

We \underline{claim} that the element $\alpha \in \cc^\times$
determined by the equation $u_{\sigma^{(0)}} u_{\pi^{(0)}}
u_{\sigma^{(0)}}^{-1} = \alpha u_{\pi^{(0)}}^3$ is a root of unity
of order $4$. To see this note that $(u_g, u_{\pi^{(0)}}^2) = 1$ for
any element $g$ in $G$ of order $2$ that commutes with $\pi^{(0)}$.
It follows that $(u_{\sigma^{(0)}}, u_{\pi^{(0)}}^2) = -1$, for
otherwise $u_{\pi^{(0)}}^2$ is contained in the center of the
algebra. This proves the claim. Therefore, we may set $\pi^{(1)} =
{\pi^{(0)}}^{3}$ if necessary and we obtain a generating set
$\Phi^{(1)} = \{\pi^{(1)},\sigma^{(1)}, \tau^{(1)} \gamma^{(1)}_1,
\dots, \gamma^{(1)}_{2m} \}$ of $G$, equivalent to $\Phi$, such that
$u_{\sigma^{(1)}} u_{\pi^{(1)}} u_{\sigma^{(1)}}^{-1} =
\varepsilon u_{\pi^{(1)}}^3$. 
By the nondegeneracy of $c$ there is an element $h^{(1)} \in
\{\tau^{(1)}, \gamma^{(1)}_1, \dots, \gamma^{(1)}_{2m}\}$ (these
elements generate $C_G(\pi^{(1)})$, the centralizer of $\pi^{(1)}$
in $G$) such that $(u_{h^{(1)}}, u_{\pi^{(1)}}) = -1$. We set
$\tau^{(2)} = h^{(1)}$ and $h^{(2)} = \tau^{(1)}$ and leave the
other generators unchanged. Next, arguing as in the previous case we
may assume that $(u_{\tau^{(2)}}, u_{\sigma^{(2)}})= 1$. Thus we
have a generating set $\Phi^{(2)}$ such that $u_{\pi^{(2)}},
u_{\tau^{(2)}}$ and $u_{\sigma^{(2)}}$ satisfy equation
(\ref{422_normal_form}). Finally, as in the previous case we obtain
a generating set $\Phi^{(r)}=\{\pi^{(r)},\sigma^{(r)}, \tau^{(r)},
\gamma^{(r)}_1, \dots, \gamma^{(r)}_{2m}\}$ (some $r$), whose
elements satisfy (a) $u_{\gamma^{(r)}_i}$ centralize the subalgebra
$\cc(u_{\pi^{(r)}}, u_{\sigma^{(r)}}, u_{\tau^{(r)}})$ for all $i$
(b) $(u_{\gamma^{(r)}_{2k-1}},u_{\gamma^{(r)}_{2k}})= - 1$ for all
$1\leq k\leq m$, and (c) all other commutators of the
$u_{\gamma^{(r)}_i}$'s are trivial. Therefore $c$ is of standard
form with respect to $\Phi^{(r)}$. This completes the proof. \qed

\begin{cor}
\label{independence.cor} If $G$ is a group on the list $\Lambda$ and
$c$ and $c'$ are nondegenerate cocycles on $G$ with values in
$\mathbb{C}^\times$, then the universal central simple algebras
$Q(U_{G,c})$ and $Q(U_{G,c'})$ are isomorphic.
\end{cor}

\proof By Theorem \ref{transitive.thm} there is an automorphism
$\varphi$ of $G$ such that $\varphi(c)$ and $c'$ are cohomologous
two-cocycles. In particular it follows that $Q(U_{G,c})$ and
$Q(U_{G,c'})$ have the same field of definition $\qq(\mu)$ and
$\varphi$ induces an automorphism $\varphi_{*}$ of $\qq(\mu)$.
Clearly, $\varphi_{*}$ extends to a $\mathbb{Q}$ isomorphism of the
free algebras $\Sigma_{c}(\mathbb{Q}(\mu))\longrightarrow
\Sigma_{\varphi(c)}(\mathbb{Q}(\mu))$ by putting $x_{ig}\longmapsto
x_{i\varphi(g)}$. Furthermore it induces an isomorphism of the
corresponding universal algebras $U_{G,c}=
\Sigma_{c}(\mathbb{Q}(\mu))/T_{c}(\mathbb{Q}(\mu))$ and $U_{G,
\varphi(c)}$, and therefore an isomorphism of their central
localization, namely the universal central simple $G$--graded
algebras.\qed

\begin{cor}
\label{divisionalg.cor}
Let $G$ be  a group of central type and let $c$ be a nondegenerate cocycle
on $G$.  The universal central simple algebra $Q(U_{G,c})$
is a division algebra if and only if $G$ is on the list $\Lambda$.
\end{cor}

\proof This follows from the the previous corollary and Proposition
\ref{list.prop}. \qed

Recall from the introduction that for a group $G$ of central type we
let $\ind(G)$ denote the maximum over all nondegenerate cocycles $c$
of the indices of the simple algebras $Q(U_{G,c})$. We have seen
that if $G$ is not on the list,  then  the universal central simple
algebra $Q(U_G)$ is not a division algebra and so the index of $G$
is strictly less than $\ord(G)^{1/2}$. In fact it is shown in
Aljadeff and Natapov \cite{AN} that the groups on the list are the
only groups responsible for the index of $Q(U_G)$. Here is the
precise result.

\begin{thm}\label{responsible.thm}
Let $P$ be a $p$-group of central type and let $Q(U_{P,c})$ be the
universal central simple $P$--graded algebra for some nondegenerate
two-cocycle $c$ on $P$. Then there is a sub-quotient group $H$ of
$P$ on the list $\Lambda$, such that $Q(U_{P,c})\cong
M_{p^{r}}(Q(U_{H,\alpha}))$, where the cocycle $\alpha$ on $H$ is
obtained in a natural way from $c$, namely there is a subgroup
$\hat{H}$ of $P$ such that $\hat{H}/N\cong H$ and
$\inf([\alpha])=\res([c])$ on $\hat{H}$.
\end{thm}

\begin{cor} Let $G$ be a central type $p$-group. Consider all
central type groups $H$ that are isomorphic to sub-quotients of $G$
and are members of $\Lambda$. Then $\ind(G)\leq
\sup(\ord(H)^{1/2})$.
\end{cor}

DeMeyer and Janusz prove in \cite[Corollary 4]{DJ} that for an
arbitrary group of central type $G$ and a nondegenerate two-cocycle
$c \in Z^2(G, \cc^\times)$ any Sylow-$p$ subgroup $G_p$ of $G$ is of
central type, and furthermore the restriction of $c$ to $G_p$ is
nondegenerate. Using this result Geva in his thesis proved that for
any subfield $k$ of $\cc$ the twisted group algebra $k^c G$ is
(non-graded) isomorphic to $\bigotimes_{p | \ord(G)} k^c G_p$. It
follows that for each prime $p$ dividing the order of $G$ there is a
$p$-group $H_{p}$ on the list $\Lambda$ which is isomorphic to a
sub-quotient of $G$ and such that $F^{s}G\cong \bigotimes_{p}
M_{p^{r}}(F^{\alpha_{p}}H_{p})$. Combining this with the preceding
Corollary we obtain:

\begin{cor}
Let $G$ be a group of central type, and for each prime $p$, $G_p$ be
a Sylow-$p$ subgroup of $G$. Consider all central type groups
$H_{p}$ that are isomorphic to sub-quotients of $G_p$ and are
members of $\Lambda_p$. Then $\ind(G)\leq \prod_{p}
\sup(\ord(H_{p})^{1/2})$.
\end{cor}

It follows from Proposition \ref{list.prop}  that for  any group $G$
of central type not on the list and any $G$--grading on
$M_n(\mathbb{C})$, the universal central simple algebra is not a
division algebra.  That means we can find nonidentity polynomials
$p(x_{ig})$ and $q(x_{ig})$ over the field of definition
$\mathbb{Q}(\mu)$ such that $p(x_{ig})q(x_{ig})$ is a graded
identity. In fact there will exist a nonidentity polynomial
$p(x_{ig})$ over  $\mathbb{Q}(\mu)$ such that $p(x_{ig})^{r}$ is a
graded identity for some positive integer $r$. This is clearly
equivalent to saying that under every substitution the value of $p$
in $M_n(\mathbb{C})$ is a nilpotent matrix.  We will refer to such a
polynomial as a {\it nilpotent} polynomial.

We present an explicit example of a nilpotent polynomial. We let
$S_3$ be the permutation group on three letters, and $C_6 = \langle
z \rangle$ be a cyclic group of order $6$. Let $\sigma = (123)$ and
$\tau = (12)$ be generators of $S_3$ and define an action of $S_3$
on $C_6$ by $\tau(z) = z^{-1}$ and $\sigma(z) = z$. We consider the
group $G = S_3 \ltimes C_6$.  Note that $G$ is not nilpotent and
hence not on the list $\Lambda$.

Let $C_6^\vee = \Hom(C_6 ,\cc^\times) = \langle \chi_z \rangle$
denote the dual group of $C_6$ and let $\langle \ , \ \rangle$
denote the usual pairing between $C_6$ and $C_6^\vee$ (i.e. $\langle
a , \chi \rangle = \chi(a)$ for all $a \in C_6, \ \chi \in
C_6^\vee$). It is well known (\cite[Corollary 10.10.2]{Ka}) that
$H^1(S_3, C_6^\vee) \cong H^2(S_3 \ltimes C_6, \cc^\times)$. The
isomorphism may be given by
\begin{equation}\label{thecocycle*}
\pi \mapsto c \ : \ c(h_1a_1, h_2a_2) = \langle h_2\cdot a_1,
\pi(h_2)\rangle,
\end{equation}
for all $h_1a_1, h_2a_2 \in S_3 \ltimes C_6$. Note that the
restriction of $c$ on $S_3$, as well as on $C_6$, is trivial.

Moreover, it is shown in Etingof and Gelaki \cite{EG1} that there
exists a bijective one-cocycle $\pi: S_3 \rightarrow C_6^\vee$ and,
because of the bijectivity, the corresponding two-cocycle $c$ is
nondegenerate. In particular this proves that $G$ is of central type
(see Etingof and Gelaki \cite{EG2}). One may define a bijective
$\pi$ as follows:
\[
\begin{array}{ll}
 \pi(\mathit{id}) = \chi_1 &  \pi(12) = \chi_{z}  \\
    \pi(123) = \chi_{z}^2    & \pi(23) = \chi_{z}^5 \\
    \pi(132) = \chi_{z}^4    & \pi(13) = \chi_{z}^3 \\
\end{array}
\]

Let $c$ denote the nondegenerate two-cocycle which corresponds to
$\pi$. Then the twisted group algebra $\cc^{c} G$ is isomorphic to
$M_6(\cc)$ or, equivalently, $M_6(\cc)$ is $G$--graded with the
class $[c]$. The group $G$ is not on the list $\Lambda$ and hence
the corresponding universal algebra $U_G$ is not a domain. Consider
the subalgebra $\cc(u_\sigma, u_y)$ where $y = z^2$. It is
straightforward to verify that the cocycle $c$ satisfies $c(y,
\sigma) = \omega^2$, where $\omega = e^{\frac{2\pi i}{3}}$ is a
primitive third root of unity in $\cc$, and hence $\cc(u_\sigma,
u_y)$ is a symbol algebra. It follows that the corresponding
subalgebra $F(t_\sigma u_\sigma, t_y u_y)$ of the universal central
simple algebra $Q(U_{G,c})=F^{s(c)}G$ is a symbol algebra $(a,b)_3$,
and by Aljadeff, et al \cite[Lemma 6]{AHN2} $a$ and $b$ are roots of
unity. It follows that $F(t_\sigma u_\sigma, t_y u_y)$ is split,
that is isomorphic to $M_3(F)$. Thus we know that there exists a
degree three nilpotent element in $U_{G,c}$. To construct it we
consider the element $u_\sigma - u_y u_\sigma$ in $\cc^{c} G$. Note
that this is not zero since $u_\sigma$ and $u_y u_\sigma$ are
linearly independent over $\cc$. Its square
\[
(u_\sigma - u_y u_\sigma)^2 = u_\sigma^2 - u_\sigma u_y u_\sigma -
u_y u_\sigma^2 + (u_y u_\sigma)^2
\]
is not zero as well, but using $u_y u_\sigma = \omega^2 u_\sigma
u_y$ and $u_y^3 = 1$ one shows that its third power
\[
(u_\sigma - u_y u_\sigma)^3 = u_\sigma^3 - u_\sigma^2 u_y u_\sigma -
u_\sigma u_y u_\sigma^2 + u_\sigma^2 u_y u_\sigma u_y u_\sigma - u_y
u_\sigma^3 + u_y u_\sigma^2 u_y u_\sigma + u_y u_\sigma u_y
u_\sigma^2 - (u_y u_\sigma)^3
\]
vanishes.

Now, as mentioned above, $res^G_{S_3}(c) = res^G_{C_6}(c) = 1$,
hence we may assume that $u_gu_h = u_{gh}$ for all $g,h \in S_3$ or
$g,h \in C_6$. In particular, we have $u_e = u_\tau^2 u_y^3$. Also,
using the cocycle values as defined in (\ref{thecocycle*}), we
obtain $u_y = \omega u_\tau u_y u_\tau u_y^2$.

Let $x_\tau = x_{1\tau}, x_\sigma = x_{1\sigma}, x_y = x_{1y}$ be
indeterminates in $\Omega=\{x_{ig} : i\in \mathbb{N}, g\in G\}$. Let
the polynomial $f(x_\tau, x_\sigma, x_y) \in
\qq(\omega)\langle\Omega\rangle$ be given by
\[
f(x_\tau, x_\sigma, x_y) = x_\sigma x_\tau^2 x_y^3 - \omega x_\tau
x_y x_\tau x_y^2 x_\sigma.
\]
It follows that the polynomial $f^3(x_\tau, x_\sigma, x_y)$ is a
polynomial identity for $M_6(\cc)$ while $f$ and $f^2$ are not.

We end this example by expressing $f^3$ explicitly as a
$\qq(\omega)-$linear combination of elementary identities:
\begin{equation*}
\begin{array}{lllll}
f^3(x_\tau, x_\sigma, x_y) & = & (x_\sigma x_\tau^2 x_y^3)^3 - (x_\tau x_y x_\tau x_y^2 x_\sigma)^3 \\
    &  +  & \omega^2 x_\tau x_y x_\tau x_y^2 x_\sigma (x_\sigma x_\tau^2 x_y^3)^2
         - \omega (x_\sigma x_\tau^2 x_y^3)^2 x_\tau x_y x_\tau x_y^2 \\
    &  +  & x_\tau x_y x_\tau x_y^2 x_\sigma (x_\sigma x_\tau^2 x_y^3)^2
         - \omega x_\sigma x_\tau^2 x_y^3 x_\tau x_y x_\tau x_y^2 x_\sigma x_\sigma x_\tau^2 x_y^3 \\
    &  +  & \omega^2 x_\sigma x_\tau^2 x_y^3 (x_\tau x_y x_\tau x_y^2 x_\sigma)^2
         - (x_\tau x_y x_\tau x_y^2 x_\sigma)^2 x_\sigma x_\tau^2 x_y^3 \\
    &  +  & \omega^2 (x_\tau x_y x_\tau x_y^2 x_\sigma) x_\sigma x_\tau^2
            x_y^3 (x_\tau x_y x_\tau x_y^2 x_\sigma)
         - \omega (x_\tau x_y x_\tau x_y^2 x_\sigma)^2 x_\sigma x_\tau^2
        x_y^3.
\end{array}
\end{equation*}

If $G$ is on the list $\Lambda$ and $c$ is a nondegenerate
two-cocycle on $G$ we have seen that the algebra $Q(U_{G,c})$ is a
division algebra of degree $n=\ord(G)^{1/2}$.  In the next result we
calculate the index of
$Q(U_{G,c})\otimes_{\mathbb{Q}(\mu)}\mathbb{C}$, where
$\mathbb{Q}(\mu)$ is the field of definition for the graded
identities.

\begin{prop}
\label{indexofext.prop} If $G$ is on the list $\Lambda$ and $c$ is a
nondegenerate two-cocycle on $G$ then the index of
$Q(U_{G,c})\otimes_{\mathbb{Q}(\mu)}\mathbb{C}$ is
$n/\ord(G^\prime)$. In particular
$Q(U_{G,c})\otimes_{\mathbb{Q}(\mu)}\mathbb{C}$ is a division
algebra if and only if $G$ is abelian of symmetric type.
\end{prop}
\proof By Aljadeff and Haile  \cite{AH} (see the discussion at the
beginning of section two) the subalgebra $F^sG^\prime$ of
$Q(U_{G,c})=F^sG$ is a cyclotomic extension of $F$.   It follows
that the index of $Q(U_{G,c})\otimes_{\mathbb{Q}(\mu)}\mathbb{C}$ is
at most $n/\ord(G^\prime)$. We proceed to prove the opposite
inequality. By Natapov \cite[proof of Theorem 6]{N} there is a
subfield $K$ of $\mathbb{C}$ containing $\mathbb{Q}(\mu)$ and a
nondegenerate cocycle $\beta$ on $G$ with values in $K$ such that
the algebra $K^\beta G\otimes_{\mathbb{Q}(\mu)}\mathbb{C}$ has index
exactly $n/\ord(G^\prime)$. By a specialization argument it follows
that the algebra $Q(U_{G,\beta})$ has index at least
$n/\ord(G^\prime)$. By Corollary \ref{independence.cor}
$Q(U_{G,\beta})$ is isomorphic to  $Q(U_{G,c})$ and hence the index
of $Q(U_{G,c})$ is at least  $n/\ord(G^\prime)$.  \qed

\end{document}